\documentclass[12pt,reqno]{amsart}
\usepackage[numbers,sort&compress]{natbib}
\usepackage{amsfonts,bbm}
\usepackage{amssymb,color}
\usepackage{amssymb}
\usepackage{amsmath,amssymb}
\usepackage{fancyhdr}
\usepackage[titletoc]{appendix}
\usepackage{enumitem}
\usepackage[colorlinks=true,backref]{hyperref}
\hypersetup{urlcolor=blue, citecolor=blue, linkcolor=black}

\usepackage[margin=3cm, a4paper]{geometry}

\usepackage{amsfonts}

\usepackage{bbm}

\usepackage{amsmath,amssymb}

\usepackage{mathrsfs}

\usepackage{hyperref}

\usepackage{graphicx}

\usepackage{float}

\usepackage{amsmath,amscd,amsbsy,amssymb,latexsym,url,bm,amsthm}

\usepackage{mathrsfs}

\usepackage{hyperref}

\usepackage{graphicx}
\usepackage{epstopdf}
\usepackage{float}

\usepackage{color}

\usepackage{amsmath,amscd,amsbsy,amssymb,latexsym,url,bm,amsthm}

\usepackage[vlined,boxed,commentsnumbered,ruled]{algorithm2e}

\usepackage{epsfig,graphicx,subfigure}

\usepackage{enumitem,balance,mathtools}

\usepackage{wrapfig}

\usepackage{mathrsfs, euscript}

\usepackage[usenames]{xcolor}

\usepackage{caption}

\usepackage{lipsum}

\usepackage{siunitx}		

\usepackage{bibentry}

\usepackage{tikz} 		

\usepackage{verbatim}	

\newtheorem{thm}{Theorem}[section]

\newtheorem{lemma}[thm]{Lemma}
\newtheorem{prop}[thm]{Proposition}
\newtheorem{defn}[thm]{ \bf{Definition}}
\theoremstyle{remark}
\newtheorem{rem}[thm]{Remark}

\makeatletter \renewenvironment{proof}[1][Proof]
{\par\pushQED{\qed}\normalfont\topsep6\p@\@plus6\p@\relax\trivlist\item[\hskip\labelsep\bfseries#1\@addpunct{.}]\ignorespaces}{\popQED\endtrivlist\@endpefalse} \makeatother
\makeatletter

\newcommand{\Cas}[1]{\begin{cases*} #1 \end{cases*}}

\newcommand{\Del}[1]{}

\newcommand{\norm}[1]{\left\|#1\right\|}		

\newcommand{\abs}[1]{\left|#1\right|}			

\newcommand{\set}[1]{\left\{#1\right\}}			

\newcommand{\inner}[1]{\left\langle#1\right\rangle} 

\newcommand{\rmnum}[1]{\uppercase\expandafter{\romannumeral#1}}  

\def\norm#1{\left\|#1\right\|}

\def\abs#1{\left|#1\right|}

\def\brk#1{\left(#1\right)}			

\def\pd{\partial}						

\def\sch{\mathscr{S}}

\def\mpq{M_{p,q}}

\def\mpqs{M_{p,q}^{s}}
\def\mpqsb{M_{p,q}^{s,\al}}

\def\eps{\varepsilon}

\def\leq{\leqslant}

\def\geq{\geqslant}

\newcommand{\N}{{\mathbb N}}			

\newcommand{\Z}{{\mathbb Z}}

\DeclareMathOperator*{\esssup}{ess\,sup}

\newcommand{\supp}{{\mbox{supp}\ }}

\def\apq{1-\rev{p}-\rev{q}}

\def\salphat{S_{\be}(t)}
\def\salpha{S_{\be}}
\newcommand{\re}{{\mathrm{Re}}}
\newcommand{\im}{{\mathrm{Im}}}

\def\al{\alpha}

\def\be{\bepha}
\def\be{\beta}

\def\th{\theta}
\def\si{\sigma}

\def\la{\lambda}
\def\ga{\gamma}
\def\real{{\mathbb{R}}}			

\def\rev#1{\frac{1}{#1}}

\def\FF{{\mathscr{F}^{-1}}}		

\def\F{{\mathscr{F}}}			

\numberwithin{equation}{section}

\begin{document}
\title[Local smoothing estimates in $\al$-modulation spaces]{Local smoothing estimates of fractional Schr\"odinger equations in $\al$-modulation spaces with some applications}

\author{Yufeng Lu}
\address{School of Mathematical Sciences\\
	Peking University\\
	Beijing 100871, China}
\email{luyufeng@pku.edu.cn}
\thanks{}



\date{}

\begin{abstract}\noindent 
We show some new local smoothing estimates of the fractional Schr\"odinger equations with initial data in $\al$-modulation spaces via decoupling inequalities. Furthermore, our necessary conditions show that the local smoothing estimates are sharp in some cases. As applications, the local smoothing estimates could show some new local well-posedness on modulation spaces of the fourth-order nonlinear Schr\"odinger equations on the line.  
\end{abstract}

\maketitle


\section{Introduction}
Denote $\salphat := \FF e^{it\abs{\xi}^{\be}} \F $, the fractional Schr\"odinger semigroup, $I=[0,1]$.
Local smoothing estimates in harmonic analysis,  based on $L^{p}$ Sobolev spaces, are as follows.
\begin{align*}
	\norm{\salphat u}_{L^{p}(I\times \real^{d})} \lesssim \norm{u}_{W^{s,p}},\ \forall u\in W^{s,p}.
\end{align*}
In contrast to the fix-time estimate: 
\begin{align*}
		\norm{S_{\be}(t_{0}) u}_{L^{p}(I\times \real^{d})} \lesssim \norm{u}_{W^{s_{0},p}},\ \forall u\in W^{s_{0},p},
\end{align*}
we know that when we take the integration of $t$ over $I$, the sharp regularity index $s$ could be less than $s_{0}$. So, we call it the local smoothing estimate. One of the significant open problems in harmonic analysis is the local smoothing conjecture of the wave equation. Until now, we only know this conjecture is true when $d=2$. One can refer to \cite{Wolff2000Local,Laba2002local,Guth2020sharp}.

In this paper, we consider the local smoothing estimate of fractional Schr\"odinger semi-group  on $\al$-modulation spaces below:\begin{align}
	\label{eq-local-smooth}
	\norm{\salphat u}_{L^{p}(I\times \real^{d})} \lesssim \norm{u}_{\mpqsb},\  \forall u\in \mpqsb.
\end{align}

The $\al$-modulation spaces $\mpqsb$, introduced by Gr\"obner in \cite{Grobner1992Banachraeume}, are proposed to be intermediate function spaces to connect modulation and Besov spaces with respect to parameters $\al\in [0,1]$. The modulation spaces $\mpqs$, introduced by Feichtinger in \cite{Feichtinger2003Modulation}, are used to measure the decay and the regularity of the function by the short time Fourier transform (STFT). Unlike the dyadic decomposition in the definition of Besov spaces, the modulation spaces can be characterized by the bounded admission partition of unity (BAPU). The $\al$-modulation spaces $\mpqsb$ could be regarded as the generalized modulation spaces. Similarly, they could be defined by $\al$ -BAPU. The precise definitions of these spaces are given in Section  \ref{sec:pre}. One can also refer to the textbooks by Gr\"ochenig \cite{Groechenig2001Foundations} and Wang et al. \cite{Wang2011Harmonic}. 

Some recent work has been devoted to the study of $\al$-modulation spaces. Borup and Nielsen in \cite{Borup2006Banach} constructed Banach frames for $\al$-modulation spaces. Kobayashi et al. in \cite{Kobayashi2009$L2$} discussed the boundedness for a class of pseudo-differential operators with symbols in $\al$-modulation spaces. Wang and Han in \cite{Han2014$$} described some basic properties of these spaces including the dual spaces, embeddings, scaling and algebraic structure. The relations between $\al$-modulation spaces and some classical function spaces such as Sobolev spaces were given by Kato in \cite{Kato2017inclusion}. One can also refer \cite{Guo2018Full,Zhao2021Sharp}. 

One of the most significant differences between $\al$-modulation spaces and Sobolev spaces is the boundedness of the unimodular Fourier multiplier such as $\salphat = \FF e^{it\abs{\xi}^{\be}} \F.$
Take for example the case of $\al=0,\be =2$. Miyachi in \cite{Miyachi1981some} showed that $S_{2}(t)$ is bounded on $L^{p}$ if and only if $p=2$. However, as showed by B\'enyi et al. in \cite{Benyi2007Unimodular}, $S_{2}(t)$ is bounded on all $\mpqs$ with $1\leq p, q \leq \infty,s\in \real$. For general cases, we refer to \cite{Miyachi2009Estimatesa,Zhao2014Remarks,Zhao2016Sharp}. In this sense, the $\al$-modulation spaces are better spaces for the initial data of the Cauchy problems for some nonlinear dispersive equations in contrast with the Sobolev spaces. One can refer to \cite{Wang2007Frequency,Guo2013Strichartz,Zhang2013Strichartz,Kato2014global,Kato2014Estimates,Kato2017applications}. Most recently, Schippa in \cite{Schippa2021smoothing} gave some sufficiency and necessity conditions of the local smoothing estimates of the Schr\"odinger semigroup in modulation spaces when $p\geq 2$ by the decoupling inequality. The goals of our paper are extending these results to the general cases with $\al, \be>0,1\leq p\leq \infty$. 


 Denote $M_{\be} = \set{(\xi,\abs{\xi}^{\be}): \xi \in \real^{d}, 1/2 \leq \abs{\xi} \leq 1}$, the compact surface in $\real^{d+1}$. One can easily know that the Gauss curvature of $M_{\be}$ is nonzero when $\be \neq 1$. The second fundamental form is positive defined when $\be >1$. 

For the case of $0<\be <1$, our main result is 
\begin{thm}
	\label{thm-al<1}
	Let $0<\be <1, 1\leq p,q  \leq \infty,s\in \real$. Denote $\al =1-\be/2 \in  (0,1),p_{0} =2+4/d$. Then \eqref{eq-local-smooth} holds if one of the following conditions is satisfied:
	\begin{description}
		\item[(A)] $p_{0} \leq p, 1/q \leq -(d+4)/dp +1 , s> \frac{\be d}{2} (1-1/p- 1/q) -\be/p$.
		\item[(B)] $2\leq p\leq p_{0}, q\geq p, s> \frac{\be d}{2} (1/2-1/q)$.
		\item[(C)] $ 1/q \geq -(d+4)/dp +1, p'\leq q \leq p, s> \frac{\be d}{4} (1-1/p -1/q).$
		\item[(D)] $q\leq p', q\leq p,s\geq 0.$
		\item[(E)] $p<2,q\geq p , s> \frac{\be d}{2}(1/p- 1/q),$
	\end{description}		
where $p'$ is the dual index of $p$ with $1/p'+1/p=1.$
\end{thm}

\begin{figure}
	\centering
	\includegraphics[width=0.7\linewidth]{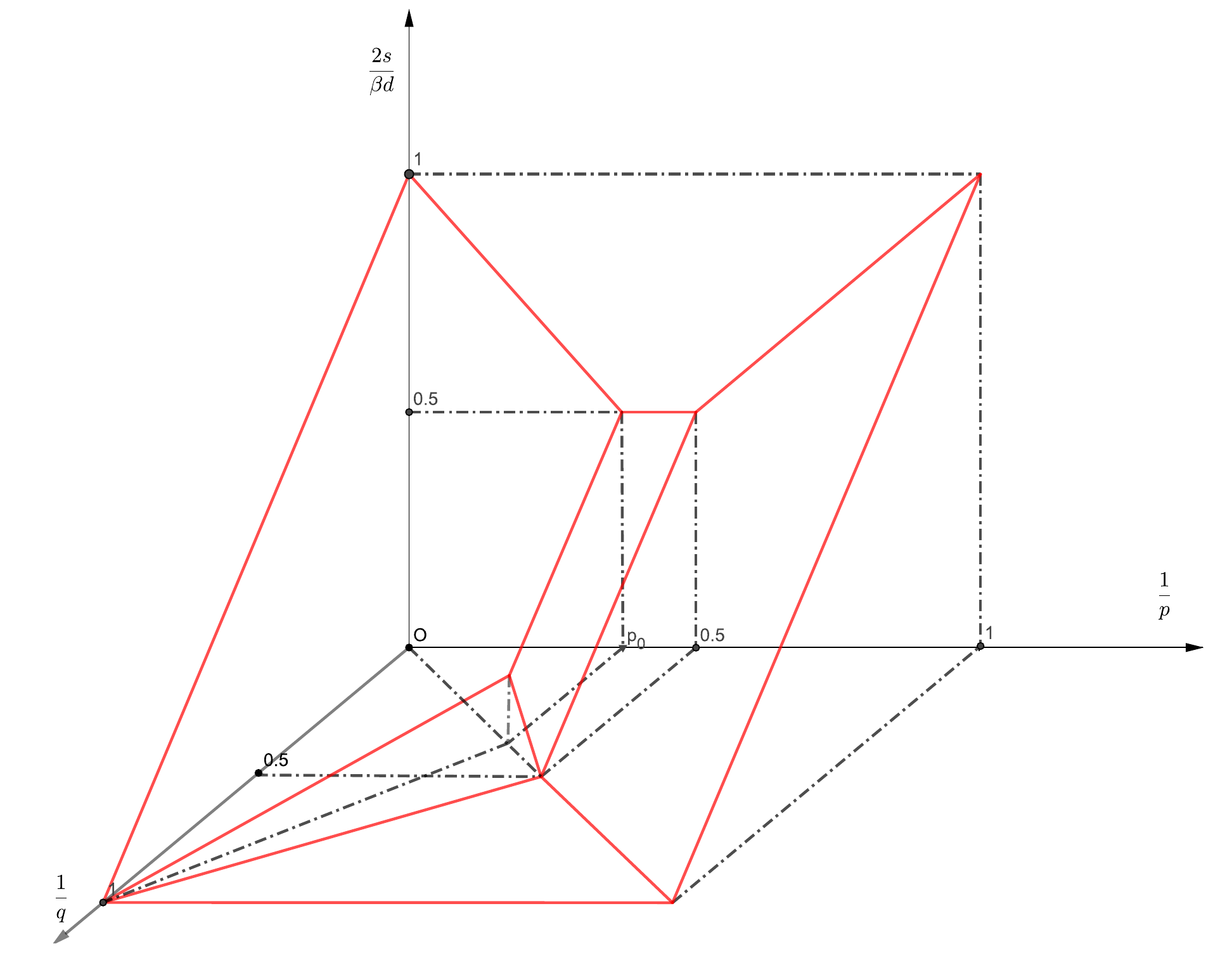}
	\caption{Relation of $(p,q,s)$ in Theorem \ref{thm-al<1}}
	\label{fig-al<1}
\end{figure}

\begin{rem}
	For the relation of $(p,q,s)$ in Theorem \ref{thm-al<1}, see Figure \ref{fig-al<1}, where $p_{0}=2+4/d$. 
\end{rem}

For the case of $1< \be \leq 2$, our main result is 

\begin{thm}
	\label{thm-1<al<2}
		Let $1<\be \leq 2, 1\leq p,q  \leq \infty,s\in \real$. Denote $\al =1-\be/2 \in  [0,1),p_{0} =2+4/d$. Then \eqref{eq-local-smooth} holds if one of the following conditions is satisfied:
		\begin{description}
			\item[(A)] $p_{0}\leq p, 1/q \leq -(d+2)/dp +1, s> \frac{\be d}{2} (1-1/p-1/q) -\be/p.$
			\item[(B)] $2\leq p\leq p_{0}, q\geq 2, s > \frac{\be d}{2} (1/2-1/q).$
			\item[(C)] $1/q \geq -(d+2)/dp +1,p'\leq q \leq 2, s>0$.
			\item[(D)] $q\leq p, q \leq p',s\geq 0$.
			\item[(E)] $p\leq 2, q>p, s>\frac{\be d}{2}(1/p- 1/q)$.
		\end{description}
\end{thm}

\begin{figure}
	\centering
	\includegraphics[width=0.7\linewidth]{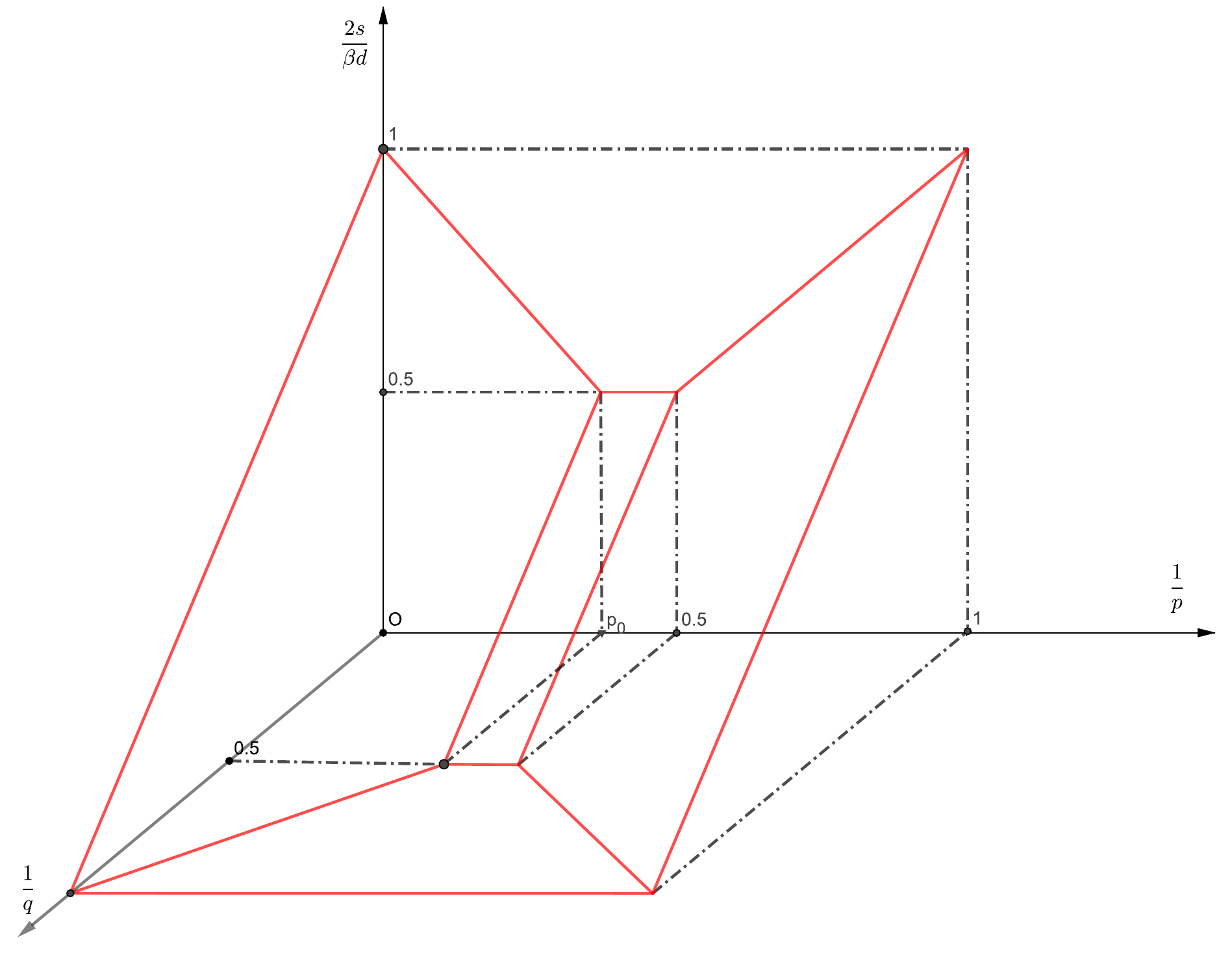}
	\caption{Relation of $(p,q,s)$ in Theorem \ref{thm-1<al<2}}
	\label{fig-al<2}
\end{figure}

\begin{rem}
	When $\be=2,\al=0$, we know that $\mpqsb = \mpqs$, the theorem above is just the Theorem 1.1 in \cite{Schippa2021smoothing}. For the relation of $(p,q,s)$ in Theorem \ref{thm-1<al<2}, see Figure \ref{fig-al<2}, where $p_{0}=2+4/d$. 
\end{rem}

\begin{rem}
	The proof of Theorem \ref{thm-1<al<2} is similar to the proof of Theorem \ref{thm-al<1}, where we use the $\ell^{2}$-decoupling instead of the $\ell^{p}$-decoupling. We omit it here for simplicity. One can also refer the proof of Theorem 1.1 in \cite{Schippa2021smoothing}.
\end{rem}

Also, we can prove the following necessary conditions when $0<\be \leq 2$.

\begin{thm}
	\label{thm-al<2-necessity}
	Let $0<\be \leq 2, 1\leq 
	p,q\leq \infty,s\in \real, $ denote $\al = 1-\be /2 \in [0,1)$. If \eqref{eq-local-smooth} holds , then we have \begin{align*}
		s\geq 0 \vee \frac{\be d}{2} (\rev{p} -\rev{q}) \vee \left(\frac{\be d}{2} (1-\rev{p} -\rev{q}) -\frac{\be}{p}\right),
	\end{align*} 
where $a\vee b = \max\set{a,b}.$
\end{thm}

\begin{rem}
	By the necessary condition above, we can see that conditions (A), (D), (E) in Theorem \ref{thm-al<1} and conditions (A),(C),(D),(E) in Theorem \ref{thm-1<al<2} are sharp without the endpoint cases.
\end{rem}

\begin{rem}
	We explain the reason of why we call the estimates in the theorems above the local smoothing estimates. Take estimates \eqref{eq-local-smooth} for example. 
	Fix $t\in I$, we consider the sharp condition of the estimate \begin{align} \label{eq-fix-time-local-smoothing}
		\norm{\salphat u}_{p} \lesssim \norm{u}_{\mpqsb} ,\ \forall u \in \mpqsb.
	\end{align}
By Theorem 4.2 in \cite{Zhao2016Sharp}, we know that $\salphat : \mpqsb \longrightarrow \mpqsb$ is bounded for any $1\leq p,q\leq \infty$ when $0<\be \leq 2, \al=1-\be/2$. Therefore, \eqref{eq-fix-time-local-smoothing} is equivalent to $\mpqsb \hookrightarrow L^{p}$. By the main result in \cite{Zhao2021Sharp} (one can also see Lemma \ref{lem-mpqsb-to-Lp} in the next section), we know that the sharp condition of this embedding is $s\geq -\be \sigma(p,q)/2$ without the endpoint case, where $\sigma(p,q)$ is denoted in the Preliminary Section below. Under this condition, take integration of $t\in I$, we can easily know that $\eqref{eq-local-smooth}$ holds. If we compare this condition with the results in Theorem \ref{thm-al<1}, we have some smoothing effects on conditions (A),(B),(C).
\end{rem}

For $\be >2$, we consider the estimates as below: 
\begin{align}
	\label{eq-local-mpq}
	\norm{\salphat u}_{L^{p}(I\times \real^{d})} \lesssim \norm{u}_{\mpqs}, \ \forall u\in \mpqs.
\end{align}
Our main result is


\begin{thm}
	\label{thm-al>2-new}
	Let $\be >2, 1\leq p,q\leq \infty,s\in \real$, denote $p_{0}=2+4/d$. Then \eqref{eq-local-mpq} holds if one of the following conditions is satisfied:
	\begin{description}
		\item[(A)] $p_{0} \leq p , 1/q \leq - (d+2)/dp +1, s> d(\be-2)(1/2-1/p) -d(1/p+1/q-1)-\be/p.$
		\item[(B)] $2\leq p \leq p_{0},q \geq 2, s>d(\be-2) (1/2 -1/p)-d(1/p+1/q-1) -\frac{\be d}{2} \brk{1/2-1/p}.$
		\item[(C)] $p_{0}\leq p,  -(d+2)/dp +1 \leq 1/q, s>(\be-2) (d/2 -(d+1)/p).$
		\item[(D)] $2\leq p \leq p_{0},q\leq 2 , s> \frac{d}{2}(\be-2) (1/2 -1/p).$
		\item[(E)] $p\leq 2,q \leq p, s\geq d(\be-2) (1/p -1/2).$
		\item[(F)] $p\leq 2, q> p, s> d(\be-2) (1/p -1/2)-d(1/q-1/p).$
	\end{description}
\end{thm}

\begin{figure}
	\centering
	\includegraphics[width=0.7\linewidth]{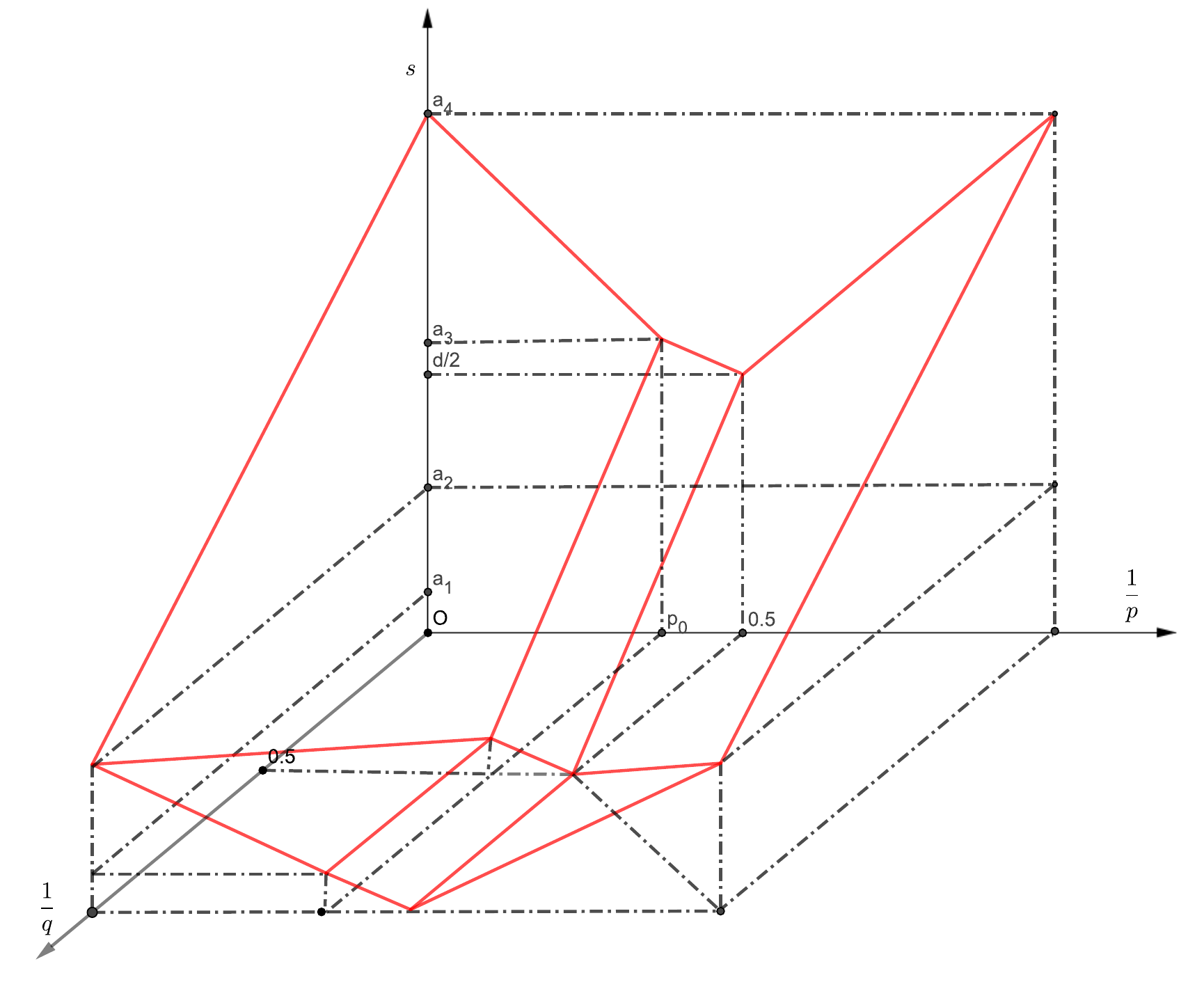}
	\caption{Relation of $(p,q,s)$ in Theorem \ref{thm-al>2-new}}
	\label{fig-al>2}
\end{figure}

\begin{rem}
	For the relation of $(p,q,s)$ in Theorem \ref{thm-al>2-new}, see Figure \ref{fig-al>2}, where $p_{0}=2+4/d, a_{1}=(\be-2)/p_{0}, a_{2}=d(\be-2)/2, a_{3}=a_{1}+d/2,a_{4}=a_{2}+d$. 
\end{rem}

\begin{rem}\label{rem-al>2-fix-time}
	As a comparison, for fix $t\in I$, we consider the sharp condition of the estimates:
	\begin{align}
		\label{eq-fix-time-al>2}
		\norm{\salphat u}_{p} \lesssim \norm{u}_{\mpqs}, \ \forall u\in \mpqs.
	\end{align}
Kobayashi and Sugimoto in \cite{Kobayashi2011inclusion} gave the almost sharp condition of this estimate. By Theorems 5.3 and 5.6 in their paper know that \eqref{eq-fix-time-al>2} holds if and only if $s\geq d(\be-2)\abs{1/2-1/p} -\sigma(p,q)$ if we omit the endpoint cases. Therefore, in Theorem \ref{thm-al>2-new}, we have some smoothing effects on conditions (A),(B),(C),(D), when we take integration of $t$ over $I$.
\end{rem}

For the necessary conditions of \eqref{eq-local-mpq}, we can get: 
\begin{thm}
	\label{thm-al>2-necessity}
	Let $\be >2,1\leq p,q \leq \infty, s\in \real$. Then if \eqref{eq-local-mpq} holds, we have $$s\geq \brk{d(1-\rev{p}-\rev{q}) -\frac{\be}{p}} \vee \brk{d(\be-2)(\rev{p}-\rev{2}) +d(\rev{p}-\rev{q})}.$$ Moreover, when $p\leq 2, p\geq  q$, we have $s\geq d(\be-2)(1/p-1/2).$
\end{thm}

\begin{rem}
	We can also consider the estimate \eqref{eq-local-mpq} when $0<\be\leq 2$, which follows by the sharp embedding between $\mpqs$ and $\mpqsb$. One can refer \cite{Guo2018Full,Han2014$$}. The results are also sharp in some cases, as in Theorems \ref{thm-al<1} and \ref{thm-1<al<2}. In fact, one can see this idea from the proof of Theorem \ref{thm-al>2-new}. For simplicity, we omit it here.
\end{rem}

The paper is organized as follows. In Section \ref{sec:pre}, we will recall some basic notation and definition, also present some useful lemmas, which will be used in our proof frequently. The proofs of our main theorems are given in Sections \ref{sec:proof of thm-al<1}-\ref{sec-al>2-nece}. Finally, we give some applications of Theorem \ref{thm-al>2-new} in partial differential equations. We get the local well-posedness of the fourth-order cubic nonlinear Schr\"odinger equations on $M_{p,2}^{s}$ for some $0<s<1/2$.

\vskip 1.5cm

\section{Preliminary}\label{sec:pre}

\vskip .5cm
	\subsection{Basic notation}

The following notation will be used throughout this article. For $0< p,q \leq \infty$, we denote 
\begin{align*}
	\sigma (p,q) &:=d\left(0 \wedge (\frac{1}{q} - \frac{1}{p}) \wedge (\frac{1}{q} +\frac{1}{p} -1)\right);\\
	\tau(p,q) & := d\left(0 \vee (\frac{1}{q} - \frac{1}{p}) \vee (\frac{1}{q} +\frac{1}{p} -1)\right),
\end{align*}
where $a\wedge b= \min\set{a,b},a\vee b =\max\set{a,b}$.

We write $\sch(\real^{d})$ to denote the Schwartz space of all complex-valued rapidly decreasing infinity differentiable functions on $\real^{d}$, and $\sch'(\real^{d})$ to denote the dual space of $\sch(\real^{d})$, all called the space of all tempered distributions. For simplification, we omit $\real^{d}$ without causing ambiguity. The Fourier transform is defined by $\F f(\xi) = \hat{f}(\xi) = \int_{\real^{d}} f(x) e^{-ix\xi} d\xi$, and the inverse Fourier transform by $\FF f(x) = (2\pi)^{-d} \int_{\real^{d}}f(\xi) e^{ix\xi} d\xi$. For $x\in \real^{d}$, we denote $\inner{x} = (1+\abs{x}^{2})^{1/2}$. For $1\leq p \leq\infty$, we define the $L^{p}$ norm:
\begin{align*}
	\norm{f}_{p} = \Cas{\left( \int_{\real^{d}} \abs{f(x)}^{p} dx \right)^{1/p}, & $1\leq p<\infty$;\\
	\esssup_{x\in \real^{d}} \abs{f(x)}, &$p= \infty$. }	
\end{align*}
We also define the $L^{p}$ Sobolev norm:
\begin{align*}
	\norm{f}_{W^{s,p}}=\norm{(I-\triangle)^{s/2}f}_{p},
\end{align*}
where $\triangle$ is the Laplace operator, $(I-\triangle)^{s/2}= \FF \inner{\xi}^{s} \F$.
Recall that the Sobolev spaces are defined by $W^{s,p}=\set{f\in \sch': \norm{f}_{W^{s,p}} <\infty}$. 

We use the notation $I \lesssim J$ if there is an independently constant C such that $I \leq C J$. Also we denote $I \approx J$ if $I \lesssim J$ and $J\lesssim I$. For $1\leq p \leq \infty$, we denote the dual index $p'$ such that $1/p+1/p'=1$.

\begin{defn}[Dyadic decomposition]
	\label{def-dyadic-decomposition}
	Choose $\psi: \real^{d} \rightarrow [0,1] $ be a smooth radial bump function adapted to the ball $B(0,2)$: $\psi(\xi) =1$ as $\abs{ \xi} \leq1$ and $\psi(\xi) =0$ as $\abs{\xi} \geq2$. We denote $\varphi(\xi) = \psi(\xi)- \psi(2\xi)$, and $\varphi_{j}(\xi) = \varphi(2^{-j} \xi)$ for $1\leq j,j\in \Z $, $\varphi_{0}(\xi ) = 1- \sum_{j\geq1} \varphi_{j}(\xi)$. Denote $\triangle_{j} = \FF  \varphi_{j} \F$.
	We say that $\set{\triangle_{j}}_{j\geq 0}$ are the dyadic decomposition operators.
\end{defn}

\vskip .5cm
\subsection{$\al$-modulation spaces}

\begin{defn}
	[$\al$-covering]\label{def-al-covering}
	Let $\al<1$. A countable set $\set{Q_{i}}_{i}$, where $Q_{i} \subseteq \real^{d}$, is called a $\al$-covering of $\real^{d}$ if:
	\begin{description}
		\item[(i)]  $\real^{d} = \cup_{i} Q_{i}$,
		\item[(ii)] $\# \set{ Q' \in Q_{i}: Q' \cap Q \neq \emptyset} \leq c(d),$ uniformly for $Q\in Q_{i}$,
		\item[(iii)] $\inner{x}^{\al d} \approx \abs{Q_{i}}$ uniformly for $x\in Q_{i}$.
	\end{description}
\end{defn}

\begin{defn}
	[$\al$-modulation spaces, \cite{Han2014$$}]\label{def alpha mpq} 
	Let $\al <1$, denote $\be=\al/(1-\al)$, suppose that $C>c>0$ are two appropriate constants such that $\set{B_{k}}_{k\in \Z^{d}}$ is a $\al$-covering of $\real^{d}$, where $B_{k} = B(\inner{k}^{\be}k, \inner{k}^{\be})$. We can choose a Schwartz function sequence $\set{\eta_{k}^{\al}}_{k\in \Z^{d}}$ satisfying 
	\begin{align*}
		\begin{cases*}
			\abs{\eta_{k}^{\al} (\xi)} \gtrsim 1,  \ \mbox{ if } \abs{\xi - \inner{k}^{\frac{\al}{1-\al}}k} < c \inner{k}^{\frac{\al}{1-\al}};\\
			\supp \eta_{k}^{\al} \subseteq \set{\xi: \abs{\xi- \inner{k}^{\frac{\al}{1-\al}}k} < C\inner{k}^{\frac{\al}{1-\al}}};\\
			\sum_{k\in\Z^{d}} \eta_{k}^{\al} (\xi) \equiv 1, \ \forall \xi \in \real^{d};\\
			\abs{\partial^{\gamma} \eta_{k}^{\al}(\xi)} \leq C_{\al} \inner{k}^{-\frac{\al\abs{\gamma}}{1-\al}}, \ \forall \xi \in \real^{d}, \gamma \in \N^{d},
		\end{cases*}
	\end{align*}
	where $C_{\al}$ is a positive constant depending only on $d$ and $\al$. We usually call these $\set{\eta_{k}^{\al}}_{k\in \Z^{d}}$ the bounded admission partition of unity corresponding ($\al-$ BAPU) to the $\al$-covering $\set{B_{k}}_{k\in \Z^{d}}$. The frequency decomposition operators can be defined by \begin{align*}
		\Box_{k}^{\al} := \FF \eta_{k}^{\al} \F.
	\end{align*}
	Let $1\leq p,q \leq \infty, s\in \real, \al\in[0,1)$, the $\al$-modulation space is defined by \begin{align*}
		M_{p,q}^{s,\al} = \set{f\in \sch': \norm{f}_{M_{p,q}^{s,\al}} = \left( \sum_{k\in\Z^{d}} \inner{k}^{sq/(1-\al)} \norm{\Box_{k}^{\al} f}_{p}^{q}  \right)^{1/q}<\infty},
	\end{align*}
	with the usual modification when $q=\infty$. 
\end{defn}

When $\al=0$, we usually denote $\mpqsb$ by $\mpqs$.

\begin{rem}
	In the previous literature, researchers usually consider $\al$-modulation spaces in the case of $\al \in [0,1)$, here we extend the definition to $\al <0$. For the convenience of readers, we put the specific definition in Appendix \ref{sec:appendix}.
\end{rem}

\subsection{Some lemmas}
In this subsection, we gather some useful results.

\begin{lemma}[\cite{Han2014$$}, embedding]
	\label{lem-embed-mpqsb}
	Let $1\leq p,q_{1},q_{2} \leq \infty, s_{1},s_{2} \in \real,\al \in [0,1)$. 
	\begin{description}
		\item[(1)] If $q_{1}>q_{2}, s_{1} + d(1-\al)/q_{1} > s_{2} + d(1-\al)/q_{2}$, then $M_{p,q_{1}}^{s_{1},\al} \hookrightarrow M_{p,q_{2}}^{s_{2},\al}$.
		\item[(2)] If $q_{1} \leq q_{2}, s_{1} \geq s_{2}$, then $M_{p,q_{1}}^{s_{1},\al} \hookrightarrow M_{p,q_{2}}^{s_{2},\al}$.
	\end{description}
\end{lemma}

\begin{lemma}[\cite{Han2014$$}, interpolation]
	\label{lem-interpolation-mpqsb}
Let $\al\in [0,1), 1\leq p,p_{0},p_{1},q,q_{0},q_{1} \leq \infty, s,s_{0},s_{1} \in \real$.	Suppose $0< \theta <1$ and 
	\begin{align*}
		s=(1-\th)s_{0} + \th s_{1},\ \frac{1}{p} = \frac{1-\th}{p_{0}} +\frac{\th}{p_{0}}, \ \rev{q} =\frac{1-\th}{q_{0}} +\frac{\th}{q_{1}},
		\end{align*}
	then we have \begin{align*}
		\brk{M_{p_{0},q_{0}}^{s_{0},\al}, M_{p_{1},q_{1}}^{s_{1},\al}}_{\th} = \mpqsb.
	\end{align*}
\end{lemma}

\begin{lemma}[\cite{Zhao2021Sharp}]
	\label{lem-mpqsb-to-Lp}
	Let $1\leq p,q\leq \infty, s\in \real,\al\in [0,1)$. Then $\mpqsb \hookrightarrow L^{p}$ if and only if one of the following conditions is satisfied:
	\begin{description}
		\item[(1)] $q\leq p < \infty, s\geq -(1-\al) \sigma(p,q);$
		\item[(2)] $p<q, s > -(1-\al) \sigma(p,q);$
		\item[(3)] $p=\infty, q=1,s \geq -(1-\al) \sigma(p,q);$
		\item[(4)] $p =\infty,q>1, s> -(1-\al) \sigma(p,q).$
	\end{description}
\end{lemma}

\begin{lemma}
	\label{lem-al>2-S(t)-FL1}
	Let $\be >2, t\in \real^{+},k\in \Z^{d}, 1\leq p\leq \infty$. Denote $\rho_{k,t}(\xi)= e^{it\abs{\xi}^{\be}} \sigma(\xi-k)$. Then we have 
	\begin{align*}
		\norm{\FF \rho_{k,t}(\cdot)}_{1} \lesssim \inner{t\abs{k}^{\be-2}}^{d/2}.
	\end{align*}
	 Then for any $u\in L^{p}, k\in \Z^{d}$, we have \begin{align*}
		\norm{\salphat \Box_{k} u}_{p} \leq C \inner{t\abs{k}^{(\be-2)}}^{d\abs{1/2-1/p}} \norm{u}_{p}		
	\end{align*}
Moreover, if we denote $\al=1-\be/2<0, \Box_{k}^{\al} = \FF \sigma(\inner{k}^{-\al/(1-\al)} \xi-k) \F$, we have \begin{align*}
	\norm{\salphat \Box_{k}^{\al} u}_{p} \lesssim \inner{t}^{d\abs{1/2-1/p}} \norm{ u}_{p}.
\end{align*}
\end{lemma}

\begin{proof}
The estimates of $\FF \rho_{k,t}$ can be found in Lemma 3.3 of \cite{Miyachi2009Estimatesa}. Notice that $\salphat \Box_{k} u = \FF \rho_{k,t} * u$, the $(p,p)$ boundedness of $\salphat \Box_{k}$ follows by the convolution of Young's inequality and the interpolation of the $(2,2)$ boundedness. As for $\salphat\Box_{k}^{\al}$, if we denote $\rho_{k,t}^{\al}(\xi) = \FF e^{it\abs{\xi}^{\be}} \sigma (\inner{k}^{-\al/(1-\al)} \xi-k)$, then $\salphat\Box_{k}^{\al}u = \FF \rho_{k,t}^{\al} * u$. Therefore, we only need to estimate $\norm{\FF \rho_{k,t}^{\al}}$. By the scaling invariance of the $\F L^{1}$ norm, we have \begin{align*}
	\norm{\FF \rho_{k,t}^{\al}}_{1} = \norm{\FF \rho_{k,\inner{k}^{\be\al/(1-\al)}t}} \lesssim \inner{\inner{k}^{\be\al/(1-\al)}t \abs{k}^{\be-2}}^{d\abs{1/2-1/p}} \lesssim \inner{t}^{d\abs{1/2-1/p}},
\end{align*}
where $ \be \al /(1-\al) +\be -2 =0$ when $\al=1-\be/2.$
\end{proof}

Taking the $\ell_{k}^{q}$ norm on both sides of the above estimates, we have 
\begin{lemma}[\cite{Miyachi2009Estimatesa}]
	\label{lem-al>2-S(t)-mpqs}
		Let $\be >2, 1\leq p,q\leq \infty, s\in \real, s\geq d(\be-2)\abs{1/2-1/p}$. Then for any $u\in \mpqs$, we have \begin{align*}
			\norm{\salphat u}_{\mpq} \lesssim \inner{t}^{d/2}\norm{u}_{\mpqs}.
		\end{align*}
\end{lemma}

\begin{lemma}[\cite{Zhao2016Sharp}]
	\label{lem-S(t)-mpqsb}
	Let $0<\be \leq 2, \be \neq 1, 1\leq p,q \leq \infty, s\in \real$, denote $ \be= 1-\be/2 \in [0,1)$. Then for any $u\in\mpqsb$, we have \begin{align*}
		\norm{\salphat u}_{\mpqsb} \lesssim \inner{t}^{d/2} \norm{u}_{\mpqsb}.
	\end{align*}
\end{lemma}

\begin{lemma}[\cite{Chaichenets2020Local}, Littlewood-Paley theory of $\mpqsb$]
	\label{lem-dyadic-of-mpq}	
	Let $\be <1, 1\leq p,q \leq \infty, s\in \real$. Then for any $u\in\mpqsb$, we have \begin{align*}
		\norm{u}_{\mpqsb} \approx \norm{2^{js} \norm{\triangle_{j} u}_{\mpq^{0,\al}}} _{\ell_{j}^{q}},
	\end{align*}
where $\triangle_{j}$ is the dyadic decomposition operator.
\end{lemma}

\begin{rem}
	In \cite{Chaichenets2020Local}, the authors presented the above argument for the modulation space, which is the special case of $\al=0$. The main observation there is that for any $k\in \Z^{d}$, the number of $j\in \N$, where $\Box_{k} \triangle_{j} \neq 0$, is bounded by a number only depends on dimension $d$. As for $\Box_{k}^{\al}$, this property is also true. One can refer Subsection 4.2 in \cite{Han2014$$}. Therefore, with the same discussion in \cite{Chaichenets2020Local}. We can get the lemma above. For simplicity, we omit the proof.
\end{rem}

\begin{lemma}[Embedding between $\al$-modulation spaces]
	\label{lem-mpqs-embed-mpqsb}
	Let $s\in \real,\al <1,0<p,q\leq \infty$. 
	\begin{description}
		\item[(1)] If $\al> 0$, then $\mpqs \hookrightarrow \mpq^{0,\al}$ if and only if $s\geq -\al \sigma(p,q)$.
		\item[(2)] If $\al <0$, then $\mpqs \hookrightarrow \mpq^{0,\al}$ if and only if $s\geq -\al \tau(p,q)$.
	\end{description}
\end{lemma}

\begin{proof}
	The proof of $(1)$  can be found in \cite{Han2014$$}, which is the special case of the embedding between $\al$-modulation spaces. For characterization of these embedding relationships, one can refer \cite{Guo2018Full}. As for $(2)$, the proof is similar to the proof of sharp embedding $\mpq^{s_{1},\al_{1}} \hookrightarrow\mpq^{s_{2},\al_{2}}$ when $\al_{1} > \al_{2}$. One can see  Appendix \ref{sec:appendix} or Section 4 in  \cite{Han2014$$} for details.
\end{proof}

\begin{lemma}
	[\cite{Bourgain2015proof,Bourgain2017Decouplings}, Decoupling]
	\label{lem-decoupling}\ 
	\begin{description}
		\item[(1)] Let $ S=\set{(\xi,\psi(\xi)): \xi \in \real^{d}} $ be a compact $C^2$ hypersurface in $\real^{d+1}$ with positive definite second fundamental form. Denote $$E_{S} f(t,x) = \int_{\real^d} e^{it\psi(\xi) +ix\xi} f(\xi) d\xi,$$
		the Fourier extension operator of $S$. Then when $p=p_{0}=2+4/d$, for any $R\geq 1$, the following estimate holds for any $\eps >0$:
		\begin{align*}
			\norm{E_{s}f(t,x)}_{L^{p}(B_{R}^{d+1})} \lesssim R^{\eps} \brk{\sum_{\Box: R^{-1/2}\mbox{-cube}} \norm{E_{S}f_{\Box}}_{L^{p}(\omega_{B_{R}^{d+1}})}^{2}}^{1/2}.
		\end{align*}
	Here $B_{R}^{d+1}$ denotes any ball in $\real^{d+1}$ with radius $R$. $\omega_{B_{R}^{d+1}}$ denotes a smooth version of the indicator function on $B_{R}^{d+1}$. 
		\item[(2)] Let $ S=\set{(\xi,\psi(\xi)): \xi \in \real^{d}} $ be a compact $C^2$ hypersurface in $\real^{d+1}$ with nonzero Gaussian curvature. The Fourier extension operator is also defined above. Then when $p=p_{0}=2+4/d$, for any $R\geq 1$, the following estimate holds for any $\eps >0$:
		\begin{align*}
			\norm{E_{s}f(t,x)}_{L^{p}(B_{R}^{d+1})} \lesssim R^{\frac{d}{2}-\frac{d+1}{p}+\eps} \brk{\sum_{\Box: R^{-1/2}\mbox{-cube}} \norm{E_{S}f_{\Box}}_{L^{p}(\omega_{B_{R}^{d+1}})}^{p}}^{1/p}.
		\end{align*}
	\end{description}
\end{lemma}

\vskip 1.5cm
\section{Proof of Theorem \ref{thm-al<1}}\label{sec:proof of thm-al<1}

We first prove the proposition below, which is the key part of the proof of Theorem \ref{thm-al<1}.

\begin{prop}
	\label{prop-al<1}
		Let $0<\be <1$, denote $\al =1-\be/2 \in  (0,1),p_{0} =2+4/d$. If $p\geq p_{0}, s >\frac{\be d}{2} (1-\frac{2}{p}) -\frac{\be}{p}$, then we have 
		\begin{align*}
			\norm{\salphat u}_{L^{p}(I\times \real^{d})} \lesssim \norm{u}_{M_{p,p}^{s,\al}}, \ \forall u \in M_{p,p}^{s,\al}. 
		\end{align*}
\end{prop}

\begin{proof}
	We only need to prove that when $s>\frac{\be d}{2} (1-\frac{2}{p}) -\frac{\be}{p}, p\geq p_{0}$, the estimate 
	\begin{align}
		\label{eq-al<1-dyadic-decomposition}
		\norm{\salphat u}_{L^{p}(I\times \real^{d})} \lesssim \lambda^{s} \norm{u}_{M_{p,p}^{0,\al}}
	\end{align}
is true for any $\la\geq 1$ and $u\in M_{p,p}^{0,\al}$ with $\supp u \subseteq \set{\xi \in \real^{d} : \la/2 \leq \abs{\xi} \leq \la}$. 

In fact, if we have the above estimates, then by the dyadic decomposition of $u$, we have $u =\sum_{j\geq 0} \triangle _{j} u$.  By triangular inequality, we have \begin{align*}
	\norm{\salphat u}_{L^{p}(I\times \real^{d})} \leq \sum_{j\geq 0} \norm{\salphat \triangle_{j} u}_{L^{p}(I\times \real^{d})}. 
\end{align*}
We can use \eqref{eq-al<1-dyadic-decomposition} to estimate $\norm{\salphat \triangle_{j} u}_{L^{p}(I\times \real^{d})}$. So we have 
\begin{align*}
	\norm{\salphat \triangle_{j} u}_{L^{p}(I\times \real^{d})} \lesssim 2^{js} \norm{\triangle_{j}u}_{M_{p,p}^{0,\al}}.
\end{align*}
Combining the two estimates above and using H'"older's inequality, we have \begin{align*}
	\norm{\salphat u}_{L^{p}(I\times \real^{d})} \lesssim \sum_{j\geq 0} 2^{js} \norm{\triangle_{j}u}_{M_{p,p}^{0,\al}} \lesssim 
	\norm{ \norm{\triangle_{j}u}_{M_{p,p}^{0,\al}} 2^{j(s+\eps)}} _{\ell_{j}^{p}},
\end{align*}
for any $\eps>0$. Here we use the H\"older's inequality. Then by Lemma \ref{lem-dyadic-of-mpq}, we know that  \begin{align*}
	\norm{ \norm{\triangle_{j}u}_{M_{p,p}^{0,\al}} 2^{j(s+\eps)}} _{\ell_{j}^{p}}\approx \norm{u}_{M_{p,p}^{s+\eps,\al}},
\end{align*}
which is the result as desired.

Then we give the proof of \eqref{eq-al<1-dyadic-decomposition}. Denote $u(x)=v_{\la}(x) = v(\la x)$, then $\widehat{u} (\xi) = \la^{-d} \widehat{v}(\la^{-1}\xi)$, which means that $\supp v \subseteq \set{\xi\in\real^{d}: 1/2\leq \abs{\xi} \leq 1}$.

By the homogeneity of $\abs{\xi}^{\be}$, it can easily be checked that $\salphat v_{\la} (x) = \salpha(\la^{\be} t) v(\la x)$. So, we have \begin{align}
	\label{eq-al<1-scaling-v}
	\norm{\salphat u}_{L^{p}(I\times \real^{d})} = \norm{\salpha(\la^{\be} t) v(\la x)}_{L^{p}(I\times \real^{d})} = \la^{-(\be+d)/p} \norm{\salphat v}_{L^{p}(\la^{\be}I\times \real^{d})}.
\end{align}
In order to estimate $\norm{\salphat v}_{L^{p}(\la^{\be}I\times \real^{d})}$ by decoupling, we need the localization of physical space $\real^{d}$. We decompose $\real^{d}$ into balls with radius $\la^{\be}$, which are boundedly overlapped. 
\begin{align*}
	\real^{d} = \bigcup B_{\la^{\be}} .
\end{align*}
So, we have \begin{align}
	\label{eq-al<1-local-in-physics}
	\norm{\salphat v}_{L^{p}(\la^{\be}I\times \real^{d})} \leq \norm{\norm{\salphat v}_{L^{p}(\la^{\be}I\times B_{\la^{\be}})}}_{\ell_{B_{\la^{\be}}}^{p}}.
\end{align}
Recall that $\salphat v (x) = E_{\be} \widehat{v} (t,x)$. By $\ell^{p}$-decoupling of $M_{\be}$ (Lemma \ref{lem-decoupling}), for any $\eps >0$, we have 
\begin{align*}
	\norm{\salphat v}_{L^{p}(\la^{\be}I\times B_{\la^{\be}})} \lesssim \la^{-\be (\frac{d+1}{p} - \frac{d}{2})+\eps} \norm{\norm{\salphat \Box_{\la^{\be/2},k}v}_{L^{p}(\omega_{\la^{\be}I}\times \omega_{B_{\la^{\be}}} )}}_{\ell_{k}^{p}},  
\end{align*}
where $\Box_{\la^{\be/2},k} u = \FF \sigma(\la^{\be/2} \xi -k )\widehat{u} (\xi)$, the decomposition of $\la^{-\be/2}$-cubes. Take this estimate into \eqref{eq-al<1-local-in-physics}, we have \begin{align}	\label{eq-al<1-decoupling}
	\eqref{eq-al<1-local-in-physics} &\lesssim \la^{-\be (\frac{d+1}{p} - \frac{d}{2})+\eps} \norm{\norm{\salphat \Box_{\la^{\be/2},k}v}_{L^{p}(\omega_{\la^{\be}I}\times \omega_{B_{\la^{\be}}} )}}_{\ell_{k,B_{\la^{\be}}}^{p}} \notag\\
	&\lesssim \la^{-\be (\frac{d+1}{p} - \frac{d}{2})+\eps} \norm{\norm{\salphat \Box_{\la^{\be/2},k}v}_{L^{p}(\omega_{\la^{\be}I}\times \real^{d} )}}_{\ell_{k}^{p}},
\end{align}
where we use the fact $\sum_{B_{\la^{\be}}} \omega_{B_{\la^{\be}}} \lesssim 1$ in the last inequality.
Recall that $\widehat{v} (\xi) = \la^{d} \widehat{u} (\la\xi)$. By rescaling, we have 
\begin{align}\label{eq-al<1-Sal-scaling}
	\salphat \Box_{\la^{\be/2},k}v (x) &= \int_{\real^d} e^{it\abs{\xi}^{\be} +ix\xi} \sigma(\la^{\be/2} \xi -k) \la^{d} \widehat{u} (\la\xi) d \xi \notag \\
	&= \int_{\real^d} e^{it\la^{-\be}\abs{\xi}^{\be} +i\la^{-1}x\xi} \sigma(\la^{\be/2-1} \xi -k)  \widehat{u} (\xi) d \xi. 
\end{align}
Notice that when $\supp \sigma(\la^{\be/2-1} \cdot -k) \cap \supp \widehat{u} \neq \emptyset$, we have $\inner{k} \approx \la^{\be/2}$. So, if we denote $\al = 1-\be/2, $ we have $\sigma(\la^{\be/2-1} \xi -k) = \sigma (\inner{k}^{-\al/(1-\al)} \xi -k)$. Taking these into \eqref{eq-al<1-decoupling} and \eqref{eq-al<1-Sal-scaling}, by scaling, we have 
\begin{align}
	\label{eq-al<1-S-alpha-mpq}
	\eqref{eq-al<1-decoupling} &\lesssim \la^{-\be (\frac{d+1}{p} - \frac{d}{2})+\eps} \norm{ \norm{\salpha(\la^{-\be} t) \Box_{k}^{\al} u (\la^{-1}x)}_{L^{p}(\omega_{\la^{\be}I} \times \real^{d})} }_{\ell_{k}^{p}} \notag\\
	&\lesssim \la^{-\be (\frac{d+1}{p} - \frac{d}{2})+\eps}  \la^{\frac{\be+d}{p}}\norm{ \norm{\salphat \Box_{k}^{\al} u}_{L^{p}(\omega_{I} \times \real^{d})} }_{\ell_{k}^{p}}.
\end{align}
By Lemma \ref{lem-S(t)-mpqsb}, we have \begin{align*}
	\norm{\salphat \Box_{k}^{\al} u}_{p} \lesssim \inner{t}^{d/2} \norm{\Box_{k}^{\al} u}_{p}.
\end{align*}
Take this into \eqref{eq-al<1-S-alpha-mpq}, we have \begin{align*}
	\eqref{eq-al<1-S-alpha-mpq} &\lesssim \la^{-\be (\frac{d+1}{p} - \frac{d}{2})+\eps}  \la^{\frac{\be+d}{p}} \norm{ \norm{ \Box_{k}^{\al} u}_{p} }_{\ell_{k}^{p}} \norm{\inner{t}^{d/2}}_{L^{p}(\omega_{I})} \\
	&\lesssim \la^{-\be (\frac{d+1}{p} - \frac{d}{2})+\eps}  \la^{\frac{\be+d}{p}} \norm{u}_{M_{p,p}^{0,\al}}.
\end{align*}
Combine the estimate above and \eqref{eq-al<1-scaling-v}-\eqref{eq-al<1-S-alpha-mpq}, we have 
\begin{align*}
		\norm{\salphat u}_{L^{p}(I\times \real^{d})} \lesssim \lambda^{-\be (\frac{d+1}{p} - \frac{d}{2})+\eps} \norm{u}_{M_{p,p}^{0,\al}},
\end{align*}
which is \eqref{eq-al<1-dyadic-decomposition} as desired.
\end{proof}

Then, we are ready to prove Theorem \ref{thm-al<1}.
\begin{proof}
	\begin{description}
		\item[(A)] By embedding of $\mpqsb$(Lemma \ref{lem-embed-mpqsb}) and Proposition \ref{prop-al<1}, when $p=p_{0},\al=1-\be/2$,  we have \begin{align*}
			\norm{\salphat u}_{L^{p}(I\times \real^{d})} \lesssim \norm{u}_{M_{p,p}^{\frac{\be d}{2} (1-\frac{2}{p}) -\frac{\be}{p}+,\al}} \lesssim \norm{u}_{M_{p,q}^{\frac{\be d}{2} (\apq) - \frac{\be}{p}+,\al}}.
		\end{align*}
	When $p=\infty,s > \frac{\be d}{2} (1-1/q) $, by Lemma \ref{lem-mpqsb-to-Lp}, we have $M_{\infty,q}^{s,\al} \hookrightarrow L^{\infty}$. Then, by Lemma \ref{lem-S(t)-mpqsb}, we have \begin{align*}
		\norm{\salphat u}_{L^{\infty} (I\times \real^{d})} \leq \norm{\salphat u}_{L_{I}^{\infty}(M_{\infty,q}^{s,\al})} \lesssim \norm{\inner{t}^{d/2}}_{L_{I}^{\infty}} \norm{u}_{M_{\infty,q}^{s,\al}} \lesssim \norm{u}_{M_{\infty,q}^{s,\al}}.
	\end{align*}
Then, the condition (A) follows by interpolation (Lemma \ref{lem-interpolation-mpqsb}) between $(p,q)=(\infty,\infty),(\infty,1),(p_{0},\infty), (p_{0},p_{0})$.
		\item[(B)] By interpolation with $p=p_{0}$, we only need to consider the case of $p=2$. Notice that \begin{align*}
			\norm{\salphat u}_{L^{2}(I\times \real^{d}) } = \norm{u}_{2}.
		\end{align*}
	So, \eqref{eq-local-smooth} is equivalent to $M_{2,q}^{s,\al} \hookrightarrow L^{2} = M_{2,2}^{0,\al}$, which is true when $s> \frac{\be d}{2} (\rev{2}-\rev{p})$ according to Lemma \ref{lem-embed-mpqsb}.
		\item[(C)] The cases of $1/q = -(d+4)/dp +1 $ and $2\leq p\leq q\leq p_{0}$ are already considered in (A)(B). By interpolation, we only need to consider the case of $1/p+1/q=1,p\geq 2$. When $s>0$, by Lemma \ref{lem-mpqsb-to-Lp}, we have $\mpqsb \hookrightarrow L^{p}$. Then,\begin{align*}
			\norm{\salphat u}_{L^{p}(I\times \real^{d})} \lesssim \norm{\salphat u}_{L^{p}_{I} (\mpqsb)} \lesssim \norm{u}_{\mpqsb}.
		\end{align*}
		\item[(D,E)] When $q\leq p, q\leq p',s>0 $ or $p<2,p<q,s> \frac{\be d}{2}(1/p -1/q)$, by Lemma \ref{lem-mpqsb-to-Lp}, we have $\mpqsb \hookrightarrow L^{p}$. Then, by the same argument in (C), we have the result as desired.
	\end{description}
\end{proof}

\vskip 1.5cm
\section{Proof of Theorem \ref{thm-al<2-necessity}} 

\begin{proof}
	Choose $\varphi \in \sch$ with $\supp \widehat{\varphi} \subseteq \set{\xi\in \real^{d}: 1/2 \leq \abs{\xi} \leq 1}$. Then we have $\norm{\varphi}_{\mpqsb} = \norm{\varphi}_{p} \approx 1$. Recall that \begin{align*}
		\salphat \varphi (x) = \int_{\real^d} e^{it\abs{\xi}^{\be} +ix\xi} \widehat{\varphi} (\xi) d\xi,
	\end{align*}
so, we know that $\abs{\salphat \varphi (x)} \gtrsim 1$ when $\abs{x} \ll 1, \abs{t} \ll 1$.

\begin{description}
	\item[(a)] For any $k\geq 10, k \in \Z^{d}$, denote $c_{k}=\inner{k}^{\rev{1-\al}}$. Take $u(x)= M_{c_{k}} \varphi (x)=e^{ic_{k} x} \varphi(x)$, then $\widehat{u} (\xi) =\widehat{\varphi}(\xi-c_{k})$ with $\supp u \subseteq c_{k} + [-1,1]^{d}$. So, we know that $\norm{u}_{\mpqsb} = \inner{k}^{\frac{s}{1-\al}}$. 
	
	Next, we estimate $\salphat u$:
	\begin{align*}
		\salphat u (x) &= \int_{\real^d} e^{it\abs{\xi}^{\be} +ix\xi} \widehat{\varphi}(\xi-c_{k}) d\xi = e^{ixc_{k}} \int_{\real^d} e^{it\abs{\xi+c_{k}}^{\be} + ix\xi} \widehat{\varphi} (\xi) d\xi\\
		&= e^{it\abs{c_{k}}^{\be} +ixc_{k}} \int_{\real^d} e^{it (\abs{\xi+c_{k}}^{\be} - \abs{c_{k}}^{\be} - \be \abs{c_{k}}^{\be-2} c_{k} \xi) + i \xi(x+\be \abs{c_{k}}^{\be-2} c_{k}t)} \widehat{\varphi} (\xi) d\xi.
	\end{align*}
Denote $h(\xi) = \abs{\xi}^{\be}$, then by Taylor's expansion, we have $\abs{\xi+c_{k}}^{\be} - \abs{c_{k}}^{\be} - \be \abs{c_{k}}^{\be-2} c_{k} \xi = h(\xi+c_{k}) -h(c_{k}) - \nabla h(c_{k}) \xi = \sum_{\abs{\ga} =2}$ $\int_{0}^{1} (1-t) \pd^{\ga} h(c_{k}+t\xi) \xi^{\ga} dt$.
Notice that when $0<\be \leq 2,$ for any $\abs{\ga}=2,\abs{\xi} \geq 1$, we have $\abs{\pd^{\ga} h(\xi)} \leq 1$. Therefore,  we have $\abs{\salphat u}  \gtrsim 1$ when $\abs{t} \ll 1, \abs{x+\be \abs{c_{k}}^{\be-2} c_{k}t} \ll 1$, which means that $\norm{\salphat u}_{L^{p}(I\times \real^{d})} \gg 1$. Combine the estimates of $u$, we have $s\geq 0 $ if \eqref{eq-local-smooth} holds.
	\item[(b)] For any $\la \gg 1$, take $u(x) = \varphi_{\la} (x) = \varphi(\la x)$, then we have $\widehat{u} (\xi) = \la^{-d} \widehat{\varphi} (\la^{-1} \xi)$, which means that $\supp \widehat{u} \subseteq \set{\xi\in \real^{d}: \la/2 \leq \abs{\xi} \leq \la}$. Denote $$\wedge _{\la} = \set{k\in\Z^{d}: \Box_{k}^{\al} u \neq 0}.$$ Then by orthogonality we know that for any $k\in \wedge_{\la}$, we have $\inner{k}^{\rev{1-\al}} \approx \la$. Also $^{\#} \wedge_{\la} \approx \la^{d(1-\al)}$. So, we have 
	\begin{align*}
		\norm{u}_{\mpqsb} &= \norm{\norm{\Box_{k}^{\al} u}_{p} \inner{k}^{s/(1-\al)}}_{\ell_{k\in \wedge_{\la}}^{q}} \\
		&\lesssim \la^{s} \norm{\norm{\varphi_{\la}}_{1} \norm{\FF \sigma(\inner{k}^{-\frac{\al}{1-\al}} \cdot -k)}_{p}}_{\ell_{k\in \wedge_{\la}}^{q}} \\
		&\lesssim \la^{s} \la^{-d} \inner{k}^{\frac{\al d}{(1-\al)p'}} \la^{(1-\al)d/q} \approx \la^{s-d+d\al/p'+d(1-\al)/q}.
	\end{align*}
By scaling, we have \begin{align*}
	\norm{\salphat u}_{L^{p}(I\times \real^{d})} &= \norm{\salpha(\la^{\be} t) \varphi(\la x)} _{L^{p}(I\times \real^{d})} \\
	&= \la^{-(\be+d)/p} \norm{\salphat\varphi}_{L^{p}(\la^{\be}I\times \real^{d})} \\
	&\gtrsim \la^{-(\be+d)/p},
\end{align*}
where we use the estimate of $\abs{\salphat \varphi}$ at the last inequality. Take the estimates of $u$ into \eqref{eq-local-smooth}, we have $s\geq d-d\al/p' -d(1-\al)/q -(\be+d)/p = \frac{\be d}{2} (\apq) -\frac{\be}{p}.$ 
	\item[(c)] Fix $N \in \N, $ which will be chosen later. For any $\la \gg 1$, denote $\wedge_{\la}= \set{k\in \Z^{d}: \inner{k}^{1/(1-\al)} \approx \la}.$  Take $u_{N} = \sum_{k\in \wedge_{\la}} T _{Nk} (M_{c_{k}} \varphi)$, where $M_{c_{k}}f(x) =f(x)e^{ixc_{k}}, T_{Nk}f(x)= f(x-Nk)$ are the modulation operator and the translation operator.  Then by the same discussion as in (a), we have the following. 
	\begin{align*}
		\norm{u_{N}}_{\mpqsb} &= \norm{\inner{k}^{s/(1-\al)} \norm{\Box_{k}^{\al} u_{N}}_{p}}_{\ell_{k\in \wedge_{\la}}^{q}} \\
		&\lesssim \la^{s} \norm{\norm{T_{Nk} (M_{c_{k}} \varphi)}_{p}}_{\ell_{k\in \wedge_{\la}}^{q}} \lesssim \la^{s} \la^{\frac{(1-\al)d}{q}}.
	\end{align*}
On the other hand, $\salphat u_{N} = \sum_{k\in \wedge_{\la}} \salphat T_{Nk} (M_{c_{k}} \varphi)$.  By the same calculation as in (a), we know that $$\abs{\salphat T_{Nk} (M_{c_{k}} \varphi)} \gtrsim 1\mbox{ in } D_{k},$$ 
where $D_{k}= \set{ (x,t): \abs{t} \ll 1, \abs{x-Nk+\be \abs{c_{k}}^{\be-2} c_{k}t} \ll1}.$
We can choose $N\gg 1$ so that $\set{D_{k}}_{k\in \wedge_{\la}}$ are disjoint. Then we have \begin{align*}
	\norm{\salphat u_{N}}_{L^{p}(I\times \real^{d})} \gtrsim \norm{\norm{\salphat T_{Nk} (M_{c_{k}} \varphi)}_{p}}_{\ell_{k\in \wedge_{\la}}^{p}} \gtrsim  \la^{\frac{(1-\al)d}{p}}.
\end{align*}
Taking the estimates of $u_{N}$ into \eqref{eq-local-smooth}, we have $$s\geq d(1-\al)(\rev{p}-\rev{q}) = \frac{\be d}{2} (\rev{p}-\rev{q}) .$$
\end{description}
\end{proof}

\vskip 1.5cm
\section{Proof of Theorem \ref{thm-al>2-new}}
We first prove the proposition below.
\begin{prop}
	\label{prop-al>2-new}
	Let $\be >2$, denote $p=p_{0} = 2+4/d$. If $s>d(\be-1)(1/2-1/p) -\be/p $, then we have \begin{align*}
		\norm{\salphat u}_{L^{p}(I\times \real^{d})} \lesssim \norm{u}_{M_{p,2}^{s}}
	\end{align*}
holds for any $u\in M_{p,2}^{s}.$
\end{prop}

\begin{proof}
	Denote $\al=1-\be/2 <0$, by the same method as in the proof of Theorem \ref{thm-al<1}, for any $\la \geq 1$  and any $u\in\sch$ with $\supp u \subseteq \set{\xi: \la/2 \leq \abs{\xi} \leq \la}$, denote $u=v_{\la}$. Then for any $\eps>0$, we have 
	\begin{align*}
		\norm{\salphat u}_{L^{p}(I\times \real^{d})} & = \la^{-\frac{\be +d}{p}} \norm{\salphat v}_{L^{p}(\la^{\be}I\times \real^{d})} \\
		&\lesssim \la^{-\frac{\be+d}{p}+\eps} \norm{\norm{\salphat \Box_{\la^{\be/2},k} v}_{L^{p}(\omega_{\la^{\be}I}\times \real^{d})} }_{\ell_{k}^{2}} \\
		&= \la^{-\frac{\be+d}{p}+\eps} \norm{\norm{\salpha(\la^{-\be}t) \Box_{k}^{\al} u (\la^{-1} x) }_{L^{p}(\omega_{\la^{\be}I}\times \real^{d})}}_{\ell_{k}^{2}}\\
		& = \la^{\eps} \norm{\norm{\salphat \Box_{k}^{\al} u}_{L^{p}(\omega_{I}\times \real^{d})}}_{\ell_{k}^{2}}.
	\end{align*}
By Lemma \ref{lem-al>2-S(t)-FL1}, for any $t\in I$, we have \begin{align*}
	\norm{\salphat \Box_{k}^{\al} u} _{p} \lesssim \inner{t}^{d(1/2-1/p)} \norm{\Box_{k}^{\al}u}_{p}.
\end{align*}
 Taking this into the estimate of $\salphat u$ above, we have \begin{align*}
	\norm{\salphat u}_{L^{p}(I\times \real^{d})} \lesssim \la ^{\eps}\norm{\norm{\Box_{k}^{\al} u}_{p}}_{\ell_{k}^{2}}=\la^{\eps} \norm{u}_{M_{p,2}^{0,\al}}.
\end{align*}
So, by dyadic decomposition, for any $u\in\sch$, we have \begin{align*}
	\norm{\salphat u}_{L^{p}(I\times \real^{d})} \lesssim \norm{u}_{M_{p,2}^{\eps,\al}}.
\end{align*}
Then by Lemma \ref{lem-mpqs-embed-mpqsb}, we know that $M_{p,2}^{s} \hookrightarrow M_{p,2}^{0,\be}$ if $s\geq -\be \tau(p,2)=(\be-1)d(1/2-1/p) -\be/p$, we can get the result as desired in the proposition.

\end{proof}
Then, we could give the proof of Theorem \ref{thm-al>2-new}.

\begin{proof}
	\begin{description}
		\item[(a)] 	By Remark \ref{rem-al>2-fix-time}, when $p\leq 2$ or $p=\infty$, we do not have any smooth effects. So, in these cases, the argument in Theorem \ref{thm-al>2-new} follows by the results in \cite{Kobayashi2011inclusion}.
		\item[(b)]  Proposition \ref{prop-al>2-new} above is just the case of $(p,q)=(p_{0},2)$. Then, by Lemma \ref{lem-embed-mpqsb}, we know that for any $\eps>0$, we have $M_{p_{0},1}^{s} \hookrightarrow M_{p_{0},2}^{s}$ and $M_{p_{0},\infty}^{s+d/2+\eps} \hookrightarrow M_{p_{0},2}^{s}$, which means that \eqref{eq-local-mpq} is true under the conditions of Theorem \ref{thm-al>2-new} when $(p,q)=(p_{0},1),(p_{0},\infty)$.
		\item[(c)] The rest of the cases follow by interpolation between (a) and (b) above. For example, Condition (A) in Theorem \ref{thm-al>2-new} can be obtained by interpolating between $(p,q) = (\infty,\infty), (\infty,1), (p_{0},2),(p_{0},\infty).$ The other conditions are the same.
	\end{description}

\end{proof}

\vskip 1.5cm
\section{Proof of Theorem \ref{thm-al>2-necessity}} \label{sec-al>2-nece}

\begin{proof}
	Choose $\varphi \in \sch$ with $\supp \widehat{\varphi} \subseteq \set{\xi\in \real^{d}: 1/2 \leq \abs{\xi} \leq 1}$. Then we have $\norm{\varphi}_{\mpqs} = \norm{\varphi}_{p} \approx 1$. Recall that \begin{align*}
	\salphat \varphi (x) = \int_{\real^d} e^{it\abs{\xi}^{\be} +ix\xi} \widehat{\varphi} (\xi) d\xi,
\end{align*}
so, we know that $\abs{\salphat \varphi (x)} \gtrsim 1$ when $\abs{x} \ll 1, \abs{t} \ll 1$.
\begin{description}
	\item[(a)] For any $\la \gg 1$, take $u(x) = \varphi_{\la} (x) = \varphi(\la x)$, then we have  $\widehat{u} (\xi) = \la^{-d} \widehat{\varphi} (\la^{-1} \xi)$, which means that $\supp \widehat{u} \subseteq \set{\xi\in \real^{d}: \la/2 \leq \abs{\xi} \leq \la}$. Denote $$\wedge _{\la} = \set{k\in\Z^{d}: \Box_{k} u \neq 0}.$$ Then by the orthogonality, we know that for any $k\in \wedge_{\la}$, we have $\inner{k}\approx \la$. Also $^{\#} \wedge_{\la} \approx \la^{d}$. So, by Young's convolution inequality, we have 
	\begin{align*}
		\norm{u}_{\mpqs} &\lesssim \la^{s} \norm{\norm{\Box_{k}u}_{p}}_{\ell_{k\in \wedge_{\la}}^{q}} \\
		&\lesssim \la^{s} \norm{\norm{u}_{1} \norm{\norm{\FF \sigma(\cdot-k)}_{p}}}_{\ell_{k\in \wedge_{\la}}^{q}} \\
		&\lesssim \la^{s-d+d/q}. 
	\end{align*}
By scaling, we have \begin{align*}
	\norm{\salphat u}_{L^{p}(I\times \real^{d})} &= \norm{\salpha(\la^{\be} t) \varphi(\la x)} _{L^{p}(I\times \real^{d})} \\&= \la^{-(\be+d)/p} \norm{\salphat\varphi}_{L^{p}(\la^{\be}I\times \real^{d})} \\
	&\gtrsim \la^{-(\be+d)/p},
\end{align*}
Take these two estimates into \eqref{eq-local-mpq}, we have $s\geq d(1-1/p-1/q) - \be/p.$
	\item[(b)] Take $u=\varphi_{\la}$ as above. By scaling, we have \begin{align} \label{eq-I(t,x)}
		\norm{\salphat u}_{L^{p}(I\times \real^{d})} &= \norm{\salpha(\la^{\be} t) \varphi(\la x)} _{L^{p}(I\times \real^{d})} \notag\\
		&=\la^{d(\be-1)/p} \norm{\salpha(\la^{\be} t) \varphi(\la^{\be} x)}_{L^{p}(I\times \real^{d})}\notag\\
		&:= \la^{d(\be-1)/p} \norm{I(t,x)}_{L^{p}(I\times \real^{d})},
	\end{align}
where $$I(t,x) = \salpha(\la^{\be} t) \varphi(\la^{\be} x) =\int_{\real^d} e^{i \la^{\be} (t\abs{\xi}^{\be} + x\xi)} \widehat{\varphi} (\xi) d\xi.$$ 
By stationary phase,  for any $t\in[1/2,1]$, we have
\begin{align*}
	\abs{I(t,x)} \lesssim \Cas{\la^{-\be d/2}, & $\abs{x} \leq 10$; \\ (\la^{\be} \abs{x})^{-100}, &$\abs{x} \geq 10.$ }
\end{align*}
So, we have $\norm{I(t,x)}_{L^{p}([1/2,1]\times \real^{d})} \lesssim \la^{-\be d/2}$ when $p=1,\infty$. By the $L^{2}$ isometric of $\salphat$, we have $\norm{I(t,x)}_{L^{2}([1/2,1]\times \real^{d})} = \la^{-\be d/2}$. Then by the convex inequality of $L^{p}$ norm: 
\begin{align*}
	\norm{I(t,x)}_{L^{2}([1/2,1]\times \real^{d})} &\leq \norm{I(t,x)}_{L^{p}([1/2,1]\times \real^{d})}^{\theta} \norm{I(t,x)}_{L^{\infty}([1/2,1]\times \real^{d})}^{1-\th},\mbox{ for } p<2, \theta=p/2 \\
		\norm{I(t,x)}_{L^{2}([1/2,1]\times \real^{d})} &\leq \norm{I(t,x)}_{L^{p}([1/2,1]\times \real^{d})}^{\theta} \norm{I(t,x)}_{L^{1}([1/2,1]\times \real^{d})}^{1-\th},\mbox{ for } p>2, \theta=p'/2 .
\end{align*}
Then we have $\norm{I(t,x)}_{L^{p}([1/2,1]\times \real^{d})} \geq \la^{-\be d/2}.$ Take this into \eqref{eq-I(t,x)},  we have \begin{align*}
	\norm{\salphat u}_{L^{p}(I\times \real^{d})} \gtrsim \la^{d(\be-1)/p-\be d/2}. 
\end{align*}
Take this estimate into \eqref{eq-local-mpq}, we have $s-d+d/q \geq d(\be-1)/p-\be d/2$, which means that $s\geq d(\be-2)(1/p-1/2) +d(1/p-1/q).$
	\item[(c)] When $p\leq 2$, for any $\abs{k} \geq 10,k\in \Z^{d}$, take $u=M_{k} \varphi$, then $\widehat{u} (\xi) = \widehat{\varphi}(\xi-k)$. So, we know $\norm{u}_{\mpqs} \approx \inner{k}^{s}.$  By change of variables, we have \begin{align*}
		\salphat u (x) = &\int_{\real^d} e^{it\abs{\xi}^{\be} + ix\xi} \widehat{\varphi}(\xi-k) d\xi \\
		&= e^{it\abs{k}^{\be} + ikx} \int_{\real^d} e^{it (h(\xi+k)-h(k)-\nabla h(k)\xi)} e^{i\xi(x+t\nabla h(k))} \varphi(\xi) d\xi,
	\end{align*}
where $h(\xi) =\abs{\xi}^{\be}.$ Denote $h_{k}(\xi) = h(\xi+k)-h(k)-\nabla h(k)\xi$, one can easily get that \begin{align}\label{eq-hk(x)}
	\abs{\pd^{\ga} h_{k}(\xi)} \lesssim \abs{k}^{\be-2}.
\end{align}

Take $L^{p}$ norm on both sides above, we have 
\begin{align}\label{eq-al>2-H(t,x)}
	\norm{\salphat u}_{L^{p}(I\times \real^{d})} = \norm{H(t,x)}_{L^{p}(I\times \real^{d})},
\end{align}
where $H(t,x) = \int_{\real^d} e^{ith_{k}(\xi) +ix\xi} \widehat{\varphi} (\xi) d\xi$. We can choose $\eta\in \sch $ such that $\eta(\xi)=1$ when $\abs{\xi} \leq 1$, and $\eta(\xi)=0,$ when $\abs{\xi} >2$. So, we have $\eta \widehat{\varphi} = \widehat{\varphi}.$  Denote $\eta_{t}(\xi) = e^{-ith_{k}(\xi)} \eta(\xi)$. So, we have $ \eta_{t}(\xi) \cdot \F H(t,\cdot) (\xi) = \widehat{\varphi}(\xi) $, then 
\begin{align}\label{eq-al>2-convolution}
	\varphi (x) = \FF \eta_{t} * H(t,\cdot) (x).
\end{align}

By the method of stationary phase and the estimates \eqref{eq-hk(x)}, we know that for any $t\in [1/2,1],$ we have \begin{align*}
	\norm{\FF \eta_{t}}_{\infty} = \norm{ \int_{\real^d} e^{i\inner{k}^{\be-2} (th_{k}(\xi)/\inner{k}^{\be-2} + x\xi) } \eta(\xi) d\xi}_{\infty} \lesssim \inner{k}^{-d(\be-2)/2}.
\end{align*}
By Young's convolution inequality, for any $t\in[1/2,1]$, we have $$\norm{\FF \eta_{t} * f}_{\infty} \lesssim \inner{k}^{-d(\be-2)/2} \norm{f}_{1}.$$
Obviously, $\norm{\FF \eta_{t} * f}_{2} = \norm{f}_{2}$.
Then, by interpolation, for any $p\leq 2$, we have $$ \norm{\FF \eta_{t} * f}_{p'} \lesssim \inner{k}^{-d(\be-2)(1/p-1/2)} \norm{f}_{p}.$$
Take this into \eqref{eq-al>2-convolution}, we have \begin{align*}
	1\approx \norm{\varphi}_{p'} \lesssim \inner{k}^{-d(\be-2)(1/p-1/2)} \norm{H(t,\cdot)}_{p}.
\end{align*}
Take integration of $t$ in $[1/2,1]$, we have $$1\lesssim \inner{k}^{-d(\be-2)(1/p-1/2)} \norm{H(t,x)}_{L^{p}([1/2,1]\times \real^{d})}.$$
Then by \eqref{eq-al>2-H(t,x)}, we have $\norm{\salphat u}_{L^{p}(I\times \real^{d})} \gtrsim \inner{k}^{d(\be-2)(1/p-1/2)}$. Then take the estimates of $u$ into \eqref{eq-local-mpq}, we have $s\geq d(\be-2)(1/p-1/2)$, which is the result as desired.
\end{description}

\end{proof}

\vskip 1.5cm
\section{Applications} \label{sec-application}
As  applications of our local smoothing estimates in Theorem \ref{thm-al>2-new}, we could get the local well-posedness of the following fourth order nonlinear Schr\"odinger equation \eqref{eq-4NLS} on modulation spaces with lower regularity request of the initial data. Also, we could get the well-posedness beyond $L^2$ spaces. Finally, by some energy argument, we could get the global well-posedness of \eqref{eq-4NLS}.

\begin{align}
	\label{eq-4NLS}\tag{4NLS}
	\Cas{i u_{t} +u_{xxxx} = -\abs{u}^{2} u, \ (t,x) \in \real^{+} \times \real,\\
	u(0,x) = u_{0},\ x\in \real.}
\end{align}

Our main results are:
\begin{thm}\label{thm-LWP}
	For any $D \in \set{M_{p,2}^{s}: 10/3 \leq p \leq 6, s>1/2-1/p}$, \eqref{eq-4NLS} is local well-posedness in $D$. In detail, for any $u_{0} \in D$, there exists $T(\norm{u_{0}}_{D}) >0$, and a unique solution $u$ of \eqref{eq-4NLS} with $u \in L_{t\in I_{T}}^{a(p)} L_{x}^{p}$, where $a(p)=8p/(3p-2),I_{T}=[0,T]$. Also, the date-to-solution map above is Lipschitz continuous.
\end{thm}

\begin{proof}
	For any $u_{0} \in D$, denote $I_{T} = [0,T], X= L_{t\in I_{T}}^{a(p)} L_{x}^{p}$, where $0<T\leq 1$, depends on $\norm{u_{0}}_{D}$, will be chosen later. Denote $ B_{X} (R) = \set{u\in X: \norm{u}_{X} \leq R}$, the closed ball in $X$, centering at $0$ with radius of $R$. Notice that when $10/3\leq p\leq 6,$ we have $a(p) \leq p$.
	
	For any decomposition of $u_{0}=u_{0}^{1} +u_{0}^{2}$, with $u_{0}^{1} \in M_{p,2}^{s}, u_{0}^{2} \in L^{2}$. By H\"older's inequality and the local smoothing estimates in Theorem \ref{thm-al>2-new}, we have 
	\begin{align}
		\label{eq-linear-Local smooth}
		\norm{S_{4}(t)u_{0}^{1}}_{X} \leq \abs{T}^{1/a(p)-1/p} \norm{ S_{4}(t)u_{0}^{1}}_{L^{p}(I\times \real^{d})} \lesssim \abs{T}^{(3p-10)/8p} \norm{u_{0}^{1}}_{M_{p,2}^{s}}.
	\end{align}
By Strichartz estimates of $S_{4}(t)$ (See Theorem 3.1 in \cite{Wang2011Harmonic}), we have 
\begin{align*}
	\norm{S_{4}(t) u_{0}^{2}}_{L^{\ga(p)}_{t} L_{x}^{p}} &\lesssim \norm{u_{0}^{2}}_{2},\\
	\norm{\mathcal{A} f}_{L^{\ga(p)}_{t} L_{x}^{p}}  &\lesssim \norm{f}_{L^{a(p)/3}_{t} L_{x}^{p/3}},
\end{align*} 
where  $(\ga(p),p)$  is the Strichartz pair, satisfying:\begin{align*}
	\frac{1}{\ga(p) } = \frac{1}{4} \brk{\rev{2}-\rev{p}}.
\end{align*}
$\mathcal{A} f$ is defined by $\mathcal{A}f(x,t) = \int_{0}^{t} S_{4}(t-\tau) f(\tau) d\tau$, is the nonlinear term. Notice that when $10/3\leq p \leq 6$, we have $a(p) \leq \ga(p)$. So, by H\"older's inequality, we have \begin{align}\label{eq-linear-strichartz}
	\norm{S_{4}(t)u_{0}^{2}}_{X} \leq \abs{T}^{1/a(p)-1/\ga(p)} \norm{ S_{4}(t)u_{0}^{2}}_{L^{\ga(p)}_{t} L_{x}^{p}} \lesssim \abs{T}^{1/4} \norm{u_{0}^{2}}_{2}.
\end{align}
Combine \eqref{eq-linear-Local smooth} and \eqref{eq-linear-strichartz}, we know that there exists $M>0$, depends only on $p,s,d$, such that 
\begin{align}\label{eq-linear-D}
	\norm{S_{4}(t) u_{0}}_{X} \leq M \norm{u_{0}}_{D}.
\end{align}
Take $f=\abs{u}^{2} u $ into the Strichartz estimate, also by H\"older's inequality, we have \begin{align}
	\label{eq-nonlinear}
	\norm{\mathcal{A} f}_{X} \leq \abs{T}^{1/4} \norm{\mathcal{A} f}_{L^{\ga(p)}_{t} L_{x}^{p}} \lesssim \abs{T}^{1/4} \norm{f}_{L^{a(p)/3}_{t} L_{x}^{p/3}} \leq C \abs{T}^{1/4} \norm{u}_{X}^{3}.
\end{align}
If we choose $T\ll 1$ such that $C\abs{T}^{1/4} (2M \norm{u_{0}}_{D})^{2} \leq 1/10$. Then for any $u_{0} \in B_{X}(2M\norm{u_{0}}_{D})$, by \eqref{eq-linear-D} and \eqref{eq-nonlinear}, we have 
\begin{align*}
	\norm{S_{4}(t)u_{0}}_{X} + \norm{\mathcal{A} \brk{\abs{u}^{2}u}} _{X} &\leq M \norm{u_{0}}_{D} + C \abs{T}^{1/4} \norm{u}_{X}^{3} \\
	&\leq M \norm{u_{0}}_{D}+1/10 \norm{u}_{X} \\
	&\leq 2M \norm{u_{0}}_{D}.
\end{align*}
Therefore, if we define $\mathscr{T} : u \rightarrow S_{4}(t) u_{0} + i \mathcal{A} \brk{\abs{u}^{2} u}$, we know that $$\mathscr{T}: B_{X}(2M\norm{u_{0}}_{D}) \longrightarrow B_{X}(2M\norm{u_{0}}_{D})$$
is contraction. Then by the contraction mapping principle, there exists $ u\in B_{X}(2M\norm{u_{0}}_{D})$, such that $u=\mathscr{T} u$, which means that $u\in X$ is the solution of \eqref{eq-4NLS}, satisfying $\norm{u}_{X} \leq 2 M \norm{u_{0}}_{D}$. So, the date-to-solution map is Lipschitz continuous.
\end{proof}

If the initial $u_{0}$ has some more regularity, out solution $u$ above is actually global. Precisely, we have 
\begin{thm}\label{thm-GWP}
	For any $E \in \set{M_{p,2}^{s}+L^2: 10/3\leq p \leq 6, s>7/2-2/p}$, if  initial data $u_{0} \in E$, then the solution $u$ in Theorem \ref{thm-LWP} can be extended to $t\in [0,\infty)$.
\end{thm}

\begin{proof}
	Our proof, named the energy argument, is based on the method in \cite{Schippa2021smoothing}, one can also see \cite{dodson2020global}.
	
	For any $u_{0} \in E$, by Theorem \ref{thm-LWP},  we have a unique solution $u\in X=L_{t\in I_{T}}^{a(p)} L_{x}^{p}$, such that $u=S_{4}(t)u_{0} +i \mathcal{A} \brk{\abs{u}^{2}u}$. Notice that when $s>7/2-2/p$, by Lemma \ref{lem-al>2-S(t)-mpqs}, we have \begin{align} \label{eq-W(t)-D}
		\norm{S_{4}(t) u_{0}} _{D} \leq \norm{S_{4}(t) u_{0}}_{M_{p,2}^{1/2-1/p+}} \lesssim \inner{T}^{1/2} \norm{u_{0}}_{M_{p,2}^{1/2-1/p+2(1/2-1/p)+}} \lesssim \inner{T}^{1/2} \norm{u_{0}}_{M_{p,2}^{s}}.
	\end{align}
As for the nonlinear term $\mathcal{A} \brk{\abs{u}^{2} u}$, by Strichartz estimates, we have \begin{align*}
	\norm{\mathcal{A} \brk{\abs{u}^{2} u}} _{L_{t\in I_{T}}^{\infty} (D)} \leq \norm{\mathcal{A} \brk{\abs{u}^{2} u}} _{L_{t\in I_{T}}^{\infty} L_{x}^{2}} \lesssim \norm{u}_{X}^{3}.
\end{align*}
So, we know that $\norm{u(T)}_{D} < \infty$.
	 By the semi-group theory, we can also solve the equation \eqref{eq-4NLS} at $t=T$. In this way, we have a maximal time $T^{*}>0$, where the solution $u$ exists in $[0,T]$ for any $T<T^{*}$. By the blow-up criteria, if $T^{*} < \infty$, it must be \begin{align*}
		\lim_{t\rightarrow T^{*}} \norm{u(t)}_{D} =\infty.
	\end{align*}
Denote $w(t) =S_{4}(t) u_{0}, v(t)= u(t)-w(t)$, by \eqref{eq-W(t)-D}, we know that $$\lim_{t\rightarrow T^{*}} \norm{W(t)}_{D} \leq C(T^{*}) < \infty.$$
Therefore, it must be \begin{align*}
	\lim_{t\rightarrow T^{*}} \norm{v(t)}_{D} = \infty.
\end{align*}
Notice that $L^{2} \hookrightarrow D$, so we have \begin{align}
	\label{eq-v(t)-blowup}
	\lim_{t\rightarrow T^{*}} \norm{v(t)}_{2} =\infty.
\end{align}
We only need to prove that \begin{align*}
	\sup_{0<t<T^{*}} \norm{v(t)}_{2} \leq C(T^{*}) < \infty.
\end{align*}
Denote $\inner{v,u}= \int_{\real} v \overline{u} dx$, the inner product in $L^{2}(\real)$.  The mass and energy of $v$ are defined below: 
\begin{align*}
	M_{v}(t) &= \frac{1}{2} \inner{v,v} ,\\
	E_{v}(t) &= \frac{1}{2} \inner{v_{xx},v_{xx}} + \rev{4} \inner{v^{2},v^{2}},\\
	\widetilde{E}_{v}(t)  &= \frac{1}{2} \inner{v_{xx},v_{xx}} + \rev{4} \inner{(v+w)^{2}, (v+w)^{2}} - \rev{4} \inner{w^{2},w^{2}}.	
\end{align*}
Note that a priori, it is not clear that the mass and energy of $v$ is finite for $t\neq 0$. The following computations are carried out for initial data from a priori class, for $u_{0} \in \sch$. This ensures all quantities to be finite and allows justifying integration by parts arguments. Since we prove bounds depending only on $\norm{u_{0}}_{M_{p,2}^{s}}$ for $p< \infty$, the arguments are a posterior justified by density of $\sch \subseteq D$ and well-posedness. 

By definition of $v=u-w$, we know that $v$ satisfies the equation below:
\begin{align}
	\label{eq-NLS_{w}} \tag{4NLSw}
	\Cas{iv_{t} + v_{xxxx} = -\abs{v+w}^{2}(v+w),\ (t,x) \in \real^{+} \times \real \\
	v(0,x) =0,\ x\in \real }
\end{align}
Then we have 
\begin{align*}
	\frac{d}{dt} M_{v}(t) = \re \inner{v_{t},v} &= \re \inner{iv_{xxxx} + i \abs{v+w}^{2}(v+x), v}\\
	&= \re\inner{iv_{xx}, v_{xx}} + \re \inner{i\abs{v+w}^{2}v,v} + \re \inner{i\abs{v+w}^{2}w,v} \\
	&=-\im  \inner{\abs{v+w}^{2}w,v} \\
	&= -\im \brk{ \inner{v\overline{v}w, v} + \inner{ vw \overline{w}, v} + \inner{\overline{v} w w,v} + \inner{w\overline{w} w,v}}.
\end{align*}
So, we have \begin{align} \label{eq-dMv(t)}
	\abs{\frac{d}{dt} M_{v}(t)} \lesssim \norm{v^{3} w}_{1} + \norm{v^{2}w^{2}}_{1} + \norm{vw^{3}}_{1}.
\end{align}
Notice that by Lemma \ref{lem-mpqsb-to-Lp}, we know that for any $p\leq 6\leq q\leq \infty$,  we have $M_{p,2}^{1/2+} \hookrightarrow M_{q,2}^{1/2+}\hookrightarrow L^{q}$, so we have 
\begin{align} \label{eq-wt-Lp}
	\norm{w(t)}_{q} = \norm{S_{4}(t) u_{0}}_{q} \lesssim \norm{ S_{4}(t) u_{0}}_{M_{p,2}^{1/2+}} \lesssim \inner{T^{*}}^{1/2} \norm{u_{0}}_{M_{p,2}^{1/2+2(1/2-1/p)+}} \lesssim _{T^{*}} 1. 
\end{align}
So, by H\"older's inequality, we have \begin{align*}
	\norm{v^{3} w} _{1} &\leq \norm{w}_{\infty} \norm{v}_{2} \norm{v}_{4}^{2} \lesssim_{T^{*}} \norm{v}_{2}^{2} + \norm{v}_{4}^{4};\\
	\norm{v^{2}w^{2}} & \leq \norm{w}_{\infty}^{2} \norm{v}_{2}^{2} \lesssim_{T^{*}} \norm{v}_{2}^{2};\\
	\norm{v w^{3}}_{1} & \leq \norm{w}_{6}^{3} \norm{v}_{2}\lesssim_{T^{*}} \norm{v}_{2} \lesssim_{T^{*}} \norm{v}_{2}^{2}+1.
\end{align*}
Take these estimates into \eqref{eq-dMv(t)}, we have \begin{align}
	\label{eq-dMv(t)-Ev}
	\abs{\frac{d}{dt} M_{v}(t)} \lesssim_{T^{*}} E_{v}(t) + M_{v}(t) +1.
\end{align}
By the definition of $\widetilde{E}_{v}(t)$, we have 
\begin{align*}
	\widetilde{E}_{v}(t)- E_{v}(t) &= \rev{4} \brk{\inner{(v+w)^{2}, (v+w)^{2}} - \inner{w^{2}, w^{2}} - \inner{v^{2}, v^{2}} } \\
	&= \inner{v^{2}, vw} + \rev{2} \inner{ v^{2}, w^{2}} + \inner{vw, vw} +\inner{vw, w^{2}}.
\end{align*}
So, by the same estimates of $w$ above, we have \begin{align*}
	\abs{\widetilde{E}_{v}(t)- E_{v}(t) } \leq \norm{v^{3} w}_{1} + \norm{v^{2} w^{2}}_{1} + \norm{v w^{3}}_{1} \lesssim_{T^{*}} \norm{v}_{2} \norm{v}_{4}^{2} + \norm{v}_{2}^{2} + \norm{v}_{2}. 
\end{align*}
Then by Cauchy-Schwartz inequality $ab\leq \eps a^{2} + b^{2}/4\eps$, we have 
\begin{align*}
	\abs{\widetilde{E}_{v}(t)- E_{v}(t) } \leq \rev{2} \norm{v}_{4}^{4} + C(T^{*}) \brk{\norm{v}_{2}^{2} +1 } \leq \rev{2} E_{v}(t) + C(T^{*}) \brk{ M_{v}(t) +1}.
\end{align*}
If we denote $A_{v}(t) = M_{v}(t) +1$, then we have 
\begin{align*}
	-\rev{2} E_{v}(t) - C(T^{*}) A_{v}(t) \leq \widetilde{E}_{v}(t)- E_{v}(t) \leq \rev{2} E_{v}(t) + C(T^{*}) A_{v}(t),
\end{align*}
which means that \begin{align*}
	\rev{2} \brk{ E_{v}(t) + 2C(T^{*}) A_{v}(t)} \leq  \widetilde{E}_{v}(t) + 2C(T^{*}) A_{v}(t) \leq \frac{3}{2} \brk{ E_{v}(t) + 2C(T^{*}) A_{v}(t)},
\end{align*}
which means that 
\begin{align}
	\label{eq-Ev-Evtilde}
	E_{v}(t) + 2C(T^{*}) A_{v}(t) \approx \widetilde{E}_{v}(t) + 2C(T^{*}) A_{v}(t) .
\end{align}
Take this into \eqref{eq-dMv(t)-Ev}, we have \begin{align}
	\label{eq-dMv(t)-Ev-tilde}
	\abs{\frac{d}{dt} M_{v}(t)} \lesssim_{T^{*}} \widetilde{E}_{v}(t) + 2C(T^{*}) A_{v}(t).
\end{align}
As for the similar estimate of $\widetilde{E}_{v}(t)$, by \eqref{eq-NLS_{w}}, we have 
\begin{align} \label{eq-dEvt-wxxxx}
	\frac{d}{dt} \widetilde{E}_{v} (t) &= \re\inner{v_{txx}, v_{xx}} + \re \inner{(v_{t}+ w_{t}) (v+w), (v+w)^{2}} -\re \inner{w_{t} w,w^{2}} \notag\\
	&= \re \inner{v_{t}, v_{xxxx}} + \re \inner{v_{t} + w_{t}, \abs{v+w}^{2} (v+w)} -\re \inner{w_{t}, \abs{w}^{2}w}\notag\\
	&= \re \inner{v_{t}, v_{xxxx} + \abs{v+w}^{2} (v+w)} + \re \inner{w_{t}, \abs{v+w}^{2} (v+w)-\abs{w}^{2}w}\notag\\
	&= \re \inner{v_{t}, iv_{t}} + \re \inner{w_{t}, \abs{v+w}^{2} (v+w)-\abs{w}^{2}w}\notag\\
	&= \re \inner{w_{t}, \abs{v+w}^{2} (v+w)-\abs{w}^{2}w}\notag\\
	&= \re \inner{ i w_{xxxx}, \abs{v+w}^{2} (v+w)-\abs{w}^{2}w}.
\end{align}
Compute that \begin{align*}
	\abs{v+w}^{2} (v+w)-\abs{w}^{2}w &= \brk{v^{2} \overline{v}} + \brk{2 v \overline{v} w + v^{2} w} + \brk{ 2vw \overline{w} + \overline{v} w^{2}} \\
	&:= E_{3,0} + E_{2,1} + E_{1,2}.
\end{align*} 
Take this into \eqref{eq-dEvt-wxxxx}, and take integration by parts, we have 
\begin{align}
	\label{eq-dEvt-E30-E21-E12}
	\frac{d}{dt} \widetilde{E}_{v} (t) &= -\im \inner{w_{xxxx}, E_{3,0} + E_{2,1} + E_{1,2}} \notag\\
	&= -\im \inner{w_{xx}, \brk{E_{3,0}}_{xx}} -\im \inner{w_{xx}, \brk{E_{2,1}}_{xx}} -\im \inner{w_{xx}, \brk{E_{1,2}}_{xx}}.
\end{align}
We estimate these three parts separately. 

For the $E_{3,0}$ part, we have 
\begin{align*}
\inner{w_{xx}, \brk{E_{3,0}}_{xx}} &= \inner{w_{xx}, \brk{2vv_{x} \overline{v} + v^{2} \overline{v}_{x}}_{x}} \\
&= \inner{w_{xx}, 2v_{x}^{2} \overline{v} + 2vv_{xx}\overline{v} + 4vv_{x} \overline{v}_{x} + v^{2} \overline{v}_{xx}}.
\end{align*}
So, by H\"older's inequality, we have \begin{align*}
	\abs{\inner{w_{xx}, \brk{E_{3,0}}_{xx}}} &\lesssim \norm{w_{xx} v_{xx} v^{2}}_{1} + \norm{w_{xx} v_{x}^{2} v}_{1}\\
	& \leq \norm{w_{xx}}_{\infty} \norm{v_{xx}}_{2} \norm{v}_{4}^{2} + \norm{w_{xx}}_{\infty} \norm{v}_{4} \norm{v_{x}}_{8/3}^{2}.
\end{align*}
By \eqref{eq-wt-Lp}, we know that for any $p\leq 6\leq q\leq \infty$, 
\begin{align*}
	\norm{w_{xx}}_{q}  \lesssim \inner{T^{*}}^{1/2} \norm{u_{0}}_{M_{p,2}^{2+1/2+2(1/2-1/p)+}} =\lesssim_{T^{*}} \norm{u_{0}}_{M_{p,2}^{7/2-2/p+}} \lesssim_{T^{*}} 1.
\end{align*}
By Gagliardo-Nirenberg inequality (See Appendix C.1 in \cite{Wang2011Harmonic}), we have \begin{align*}
	\norm{v_{x}}_{8/3} \lesssim \norm{v}_{4}^{1/2} \norm{v_{xx}}_{2}^{1/2}.
\end{align*}
Take this into the  estimate of $E_{3,0}$, we have \begin{align*}
	\abs{\inner{w_{xx}, \brk{E_{3,0}}_{xx}}} \lesssim_{T^{*}} \norm{v_{xx}}_{2} \norm{v}_{4}^{2} \lesssim_{T^{*}} \norm{v_{xx}}_{2}^{2} + \norm{v}_{4}^{4} \lesssim_{T^{*}} E_{v}(t).
\end{align*}

For the $E_{2,1}$ part, by the same calculation, we have \begin{align*}
	\abs{\inner{w_{xx}, \brk{E_{2,1}}_{xx}}} &\lesssim \norm{w_{xx}v_{xx}vw}_{1} + \norm{w_{xx} v_{x}^{2} w}_{1} + \norm{w_{xx} vv_{x}w_{x}}_{1} + \norm{w_{xx}^{2} v^{2}}_{1} \\ 
	& \lesssim_{T^{*}} \norm{v_{xx}}_{2} \norm{v}_{2} + \norm{v_{x}}_{2}^{2} + \norm{v}_{2} \norm{v_{x}}_{2} + \norm{v}_{2}^{2} \\
	&\lesssim_{T^{*}} \norm{v}_{2}^{2} + \norm{v_{xx}}_{2}^{2} + \norm{v_{x}}_{2}^{2}.
\end{align*}
Also, by Gagliardo-Nirenberg inequality, we have \begin{align*}
	\norm{v_{x}}_{2} \lesssim \norm{v}_{2} \norm{v_{xx}}_{2}.
\end{align*}
Therefore, we have \begin{align*}
	\abs{\inner{w_{xx}, \brk{E_{2,1}}_{xx}}} &\lesssim_{T^{*}} \norm{v}_{2}^{2} + \norm{v_{xx}}_{2}^{2} \lesssim_{T^{*}} E_{v}(t) + M_{v}(t).
\end{align*}

For the $E_{1,2}$ part, we have 
\begin{align*}
	\abs{\inner{w_{xx}, \brk{E_{1,2}}_{xx}}} &\lesssim \norm{w_{xx}v_{xx} w^{2}}_{1} + \norm{w_{xx}v_{x}ww_{x}}_{1} + \norm{w_{xx}vw_{x}^{2}}_{1} + \norm{w_{xx}^{2} vw}_{1} \\
	&\lesssim \norm{v_{xx}}_{2} \norm{w_{xx}}_{6} \norm{w}_{6}^{2} + \norm{v_{x}}_{2} \norm{w}_{6} \norm{w_{x}}_{6} \norm{w_{xx}}_{6} \\
	&+ \norm{v}_{2}  \norm{w_{x}}_{6}^{2} \norm{w_{xx}}_{6}
	+\norm{v}_{2} \norm{w}_{6}  \norm{w_{xx}}_{6}^{2}\\
	& \lesssim_{T^{*}} \norm{v_{xx}}_{2} + \norm{v_{x}}_{2} + \norm{v}_{2}\\
	&\lesssim_{T^{*}} \norm{v_{xx}}_{2}+\norm{v}_{2} \norm{v_{xx}}_{2}+ \norm{v}_{2}\\
	&\lesssim_{T^{*}} \norm{v_{xx}}_{2}^{2} +\norm{v}_{2}^{2} +1 \\
	&\lesssim_{T^{*}} E_{v}(t) + M_{v}(t) +1.
\end{align*}

Take all these three parts into \eqref{eq-dEvt-E30-E21-E12}, we have 
\begin{align*}
	\abs{\frac{d}{dt} \widetilde{E}_{v} (t)} \lesssim_{T^{*}} E_{v}(t) + M_{v}(t) +1.
\end{align*}
Also, by \eqref{eq-Ev-Evtilde}, we have 
\begin{align*}
	\abs{\frac{d}{dt} \widetilde{E}_{v} (t)} \lesssim_{T^{*}} \widetilde{E}_{v} (t) + 2C(T^{*}) A_{v}(t).
\end{align*}
Combine this estimate with \eqref{eq-dMv(t)-Ev-tilde}, we have 
\begin{align*}
	\frac{d}{dt} \brk{\widetilde{E}_{v} (t) + 2C(T^{*}) A_{v}(t) } &\leq C(T^{*}) \brk{\widetilde{E}_{v} (t) + 2C(T^{*}) A_{v}(t)},\\
	\widetilde{E}_{v} (0) + 2C(T^{*}) A_{v}(0) = 2C(T^{*}).
\end{align*}
Then, by Gronwall's inequality, we know that \begin{align*}
	\widetilde{E}_{v} (t) + 2C(T^{*}) A_{v}(t) \lesssim_{T^{*}} 1,
\end{align*}
which means that $M_{v}(t) \lesssim_{T^{*}} 1$, is contraction with \eqref{eq-v(t)-blowup}.
\end{proof}

\vskip 1.5cm
\section*{Acknowledgments}
Thanks to Robert Schippa for some helpful discussions. Thanks to my advisor Baoxiang Wang for inspiring advises.

\newpage
\begin{appendices}

\section{$\mpqsb$ with $\al<0$} \label{sec:appendix}
We first give the existence of the $\al$-covering of $\real^{d}$ when $\al<0$.
\begin{prop}
	For $\al <0$, there exists a $\al$-covering $\set{B_{k}}_{k\in \Z^{d}}$ of  $\real^{d}$ with a $\al$-BAPU $\set{\eta_{k}^{\al}}_{k\in \Z^{d}}$.
\end{prop}

\begin{proof}
	Our proof is based on Subsection 2.1 in \cite{Borup2006Banach}. Denote $\be=\al/(1-\al) \in (-1,0)$. We could define a homeomorphism $\delta_{\be}$ of $\real^{d}$ by $\delta_{\be} (x) = \abs{x}^{\be} x$. It is easy to check that $\delta_{\be}^{-1} = \delta_{\be'}$, where $\be'= -\be/(1+\be)$. Also, we know that $\delta_{\be} \in C^{1}(\real^{d} \setminus \set{0})$.
	
	For any $k\in \Z^{d}, \abs{k} \geq 2\sqrt{d}$, for any $x \in B(k,\sqrt{d})$, we know that $\abs{ x} \approx \abs{k} \gtrsim 1$. Then by the mean value inequality, we have \begin{align*}
		\abs{\delta_{\be}(x) - \delta_{\be}(k)} \leq C_{\be} \sup_{0\leq t \leq 1} \abs{tx+(1-t)k}^{\be} \abs{x-k} \leq C_{\be} \inner{k}^{\be}.
	\end{align*}
	So, we know that $\delta_{\be} \brk{B (k,\sqrt{d})} \subseteq B(\delta_{\be} (k), C_{\be} \inner{k}^{\be}) , \forall \abs{k} \geq 2\sqrt{d}$.  
	
	When $\abs{k} \leq 2 \sqrt{d}$, it is obvious that $\delta_{\be} \brk{B (k,\sqrt{d})}$ is bounded. So, we could choose $C_{\be}'>0$, such that $\delta_{\be} \brk{B (k,\sqrt{d})} \subseteq B(\delta_{\be} (k), C_{\be} \inner{k}^{\be}) , \forall \abs{k} \leq  2\sqrt{d}$. 
	
	Combine all the cases, we know that there exists $r>0$, such that 
	\begin{align*}
		\delta_{\be} \brk{B (k,\sqrt{d})} \subseteq B(\delta_{\be} (k), r \inner{k}^{\be}) , \forall k\in \Z^{d}.
	\end{align*}
	Denote $B_{k} = B(\delta_{\be} (k), r \inner{k}^{\be})$, with the condition above, we know that $\real^{d} = \bigcup_{k} B_{k}$ because that $\set{B (k,\sqrt{d})}_{k\in \Z^{d}}$ is a covering of $\real^{d}$ and $\delta_{\be}$ is a homeomorphism.
	
	Similarly, for the $r>0$ above, we have $R>0$, such that 
	\begin{align*}
		\delta_{\be'} \brk{B(\delta_{\be} (k), r \inner{k}^{\be})} \subseteq B (k,R) , \forall k\in \Z^{d},
	\end{align*}
	which is equivalent to \begin{align*}
		B(\delta_{\be} (k), r \inner{k}^{\be}) \subseteq \delta_{\be}\brk{ B (k,R)} , \forall k\in \Z^{d}.
	\end{align*}
	Then by the bounded overlap of $\set{B(k,R)}_{k\in \Z^{d}}$, we know that $\set{B_{k}}_{k\in \Z^{d}}$ is also bounded overlapped. 
	
	One can easily check that condition (iii) in Definition \ref{def-al-covering} holds for $B_{k}$ given above. So, we construction a $\al$-covering of $\real^{d}$, named $\set{B_{k}}_{k\in \Z^{d}}$. The corresponding $\al$-BAPU of this $\al$-covering can be given as follows.
	
	For any $\abs{k} \neq 0$, we can choose $r_{0} >0$, such that $r\inner{k}^{\be} \leq r_{0} \abs{k}^{\be}$. Take $\rho \in \sch$ with $\rho(\xi)=1$, when $\abs{\xi} \leq r_{0}$, $\rho(\xi) = 1$ when $\abs{\xi} \geq 2r_{0}$. Denote $\rho_{0}=\rho, \rho_{k} (\xi)= \rho(\abs{k}^{-\be} \xi -k)$ for any $\abs{k} \neq 0$. So, we know that $\rho_{k}(\xi) =1 $ when $\xi \in B_{k}$ and $\supp \rho_{k} \subseteq B(\delta_{\be}(k), 2r_{0} \inner{k}^{\be}) := B_{k}'$. Notice that $\set{B_{k}'}_{k\in \Z^{d}}$ is  bounded overlapped and $\set{B_{k}}_{k\in \Z^{d}}$ is a covering of $\real^{d}$. We have \begin{align*}
		1\leq \sum_{k\in\Z^{d}} \rho_{k} (\xi)\leq C(d).
	\end{align*}
	Denote \begin{align*}
		\psi_{k} = \frac{\rho_{k}}{\sum_{k\in\Z^{d}} \rho_{k}},
	\end{align*}
	one can check that $\set{\psi_{k}}_{k\in \Z^{d}}$ is a $\be$-BAPU.
\end{proof}

\begin{prop} \label{prop-mpqsb-embed-mpq}
	Let $\al <0, s\in \real, 0<p,q\leq \infty$. Then $\mpqs \hookrightarrow \mpq^{0,\al}$ if and only if $s\geq -\al \tau(p,q)$.
\end{prop}
\begin{proof}
	The proof is based on Section 4 in  \cite{Han2014$$}.
	
	Sufficiency:\begin{description}
		\item[(1)] When $(p,q)=(2,2), \tau(p,q) =0$, we have $L^{2} = M_{2,2} = M_{2,2}^{0,\al}$.
		\item[(2)] For any $k,\ell \in \Z^{d}$, denote \begin{align*}
			\wedge_{k} &= \set{\ell \in \Z^{d}: \Box_{\ell} \Box_{k}^{\al} \neq 0},\\
			\vee_{\ell} &= \set{k\in \Z^{d}: \Box_{\ell} \Box_{k}^{\al} \neq 0}.
		\end{align*}
		By the orthogonality of $\set{\Box_{\ell}}_{\ell\in \Z^{d}}$ and $\set{\Box_{k}^{\al}}_{k\in \Z^{d}}$, we know that \begin{align*}
			\# \wedge_{k} \lesssim_{d} 1, \ \ \# \vee_{\ell} \lesssim_{d} \inner{\ell}^{-\al d}.
		\end{align*}
		When $p=1,$ or $\infty$, we know $\tau(p,q) = d/q$. By the convolution Young's inequality and the orthogonality of these decompositions, we have 
		\begin{align*}
			\norm{u}_{\mpq^{0,\al}} &= \| \|\Box_{k}^{\al} u\|_{p}\|_{\ell_{k}^{q}}=  \norm{ \|\|\Box_{k}^{\al} u\|_{p}\|_{\ell_{k\in\vee_{\ell}}^{q}}}_{\ell_{\ell}^{q}}\\
			&= \norm{ \|\|\sum_{\ell'\in \wedge_{k} }\Box_{\ell'}\Box_{k}^{\al} u\|_{p}\|_{\ell_{k\in\vee_{\ell}}^{q}}}_{\ell_{\ell}^{q}}\\
			& \lesssim \norm{\|\sum_{\abs{\ell-\ell'} \lesssim 1}\norm{\Box_{\ell'}u}_{p}\|_{\ell_{k\in\vee_{\ell}}^{q}}}_{\ell_{\ell}^{q}}\\
			&\lesssim \norm{\norm{\inner{\ell}^{\al d/q} \Box_{\ell} u}_{p}}_{\ell_{\ell}^{q}} = \norm{u}_{\mpq^{-\al d/q}}.
		\end{align*}
		When $0<p<1$, we know $\tau(p,q) = d(1/p+1/q-1)$. By using the convolution inequality below instead of the convolution Young's inequality, we could get the embedding as well. 
		\begin{align*}
			\norm{\Box_{\ell} \Box_{k}^{\al} u} _{p} \lesssim \norm{\FF \eta_{k}^{\al}}_{p} \norm{\Box_{\ell} u}_{p} \lesssim \inner{k}^{\al d(1-1/p) /(1-\al)} \norm{\Box_{\ell} u}_{p}.
		\end{align*}
		\item[(3)] When $p=2, q<2$, we know $\tau(2,q) = d(1/q-1/2)$. By the Plancherel formula and H\"older's inequality, we have \begin{align*}
			\norm{u}_{M_{2,q}^{0,\al}} &= \| \|\Box_{k}^{\al} u\|_{2} \|_{\ell_{k}^{q}} = \norm{\norm{ \|\|\Box_{k}^{\al} u\|_{2}\|_{\ell_{k}^{q}}}_{\ell_{\ell}^{q}}}_{\ell_{\ell}^{q}} \\
			&\leq \| \|\Box_{k}^{\al} u\|_{2} \|_{\ell_{k \in \vee_{\ell}}^{q}} = \norm{ \|\|\Box_{k}^{\al} u\|_{2}\|_{\ell_{k\in \vee_{\ell}}^{2}} (\# \vee_{\ell})^{1/q-1/2}}_{\ell_{\ell}^{q}}\\
			&\lesssim \norm{\norm{\Box_{\ell} u}_{2} \inner{\ell}^{-\al d (1/q-1/2)}}_{\ell_{\ell}^{q}}=\norm{u}_{M_{2,q}^{-\al\tau(2,q)}}.
		\end{align*}
		\item[(4)] Take interpolation of (1,2,3), we could get the result for the rest cases of $(p,q)$ as desired.
	\end{description}
	
	Necessity: If we know $\mpqs \hookrightarrow \mpq^{0,\al}$, then we have 
	\begin{align} \label{eq} 
		\norm{u}_{\mpq^{0,\al}} \lesssim \norm{u}_{\mpqs}.
	\end{align}
	\begin{description}
		\item[(a)] For any $k\in \Z^{d}$, take $u= \FF \si(\inner{k}^{-\al/(1-\al)} \xi -k )$ into the inequality above.  We know that $\norm{u}_{\mpq^{0,\al}} \approx \norm{u}_{p} = \inner{k}^{\al d(1-1/p)/(1-\al)}$. Also, $\norm{u}_{\mpqs} \approx \inner{k}^{s/(1-\al)} \norm{u}_{p}.$ Then we have $s\geq 0$.
		\item[(b)] For any $\ell \in \Z^{d}$, take $u= \FF \si(\xi-\ell)$ into the inequality above. We know that $\norm{u}_{\mpqs} \approx \inner{\ell}^{s}$. Also, we have \begin{align*}
			\norm{u}_{\mpq^{0,\al}} = \norm{ \norm{\Box_{k}^{\al} u}_{p}}_{\ell_{k\in \vee_{\ell}}^{q}} \gtrsim \norm{\norm{\FF \eta_{k}^{\al}}_{p} }_{\ell_{k\in \vee_{\ell}}^{q}} \approx \inner{\ell}^{\al d(1-1/p-1/q)}.
		\end{align*}
		So, we have $s\geq -\al d(1/p+1/q-1)$.
		\item[(c)] For any $\ell \in \Z^{d}, N\in \N^{+}$, take $$u^{N} = \sum_{k \in \vee_{\ell}} T_{Nk} \brk{ \FF \si(\inner{k}^{-\al/(1-\al)} \xi -k )},$$ where $T_{Nk}f(x) = f(x-Nk)$ is the translation operator. Then we have \begin{align*}
			\norm{u^{N}}_{\mpq^{0,\al}} = \norm{  \norm{\FF  \si(\inner{k}^{-\al/(1-\al)} \xi -k )}_{p}}_{\ell_{k\in \vee_{\ell}}^{q}} \approx \inner{\ell}^{\al d(1-1/p-1/q)}.
		\end{align*}
		On the other hand, we have $\norm{u^{N}}_{\mpqs} \approx \inner{\ell}^{s} \norm{f^{N}}_{p}$. Take $N\rightarrow \infty$, by the almost orthogonality of $T_{Nk} \brk{ \FF \si(\inner{k}^{-\al/(1-\al)} \xi -k )}$, we have \begin{align*}
			\lim_{N\rightarrow \infty} \norm{u^{N}}_{p}& = \norm{ \norm{ \FF \si(\inner{k}^{-\al/(1-\al)} \xi -k )}_{p}}_{\ell_{k \in \vee_{\ell}}^{p}} \\
			&\approx \inner{\ell}^{\al d(1-1/p)} \inner{\ell}^{-\al d/p} = \inner{\ell}^{\al d(1-2/p)}.
		\end{align*} 
		Take these two estimates into the inequality \eqref{eq}, we have $s\geq -\al d(1/q-1/p)$.
	\end{description}
	Combine the three conditions above, we have $s\geq -\al \tau(p,q)$, as desired.
\end{proof}

As for $\mpqsb \hookrightarrow \mpq$, we also have 
\begin{prop}
	Let $\al <0, s\in \real,0<p,q \leq \infty$. Then $\mpqsb \hookrightarrow \mpq$ if and only if $s\geq \al \si(p,q).$
\end{prop}
\begin{proof}
	The necessity is the same as the proof of Proposition \ref{prop-mpqsb-embed-mpq}. The sufficiency can be done by the interpolation of the following cases. 
	\begin{description}
		\item[(a)] When $p=2,q>2$, we know that $\sigma(2,q)=d(1/q-1/2).$ By the Plancherel's formula and the H\"older's inequality of the sequences, we have 
		\begin{align*}
			\norm{u}_{M_{2,q}} &= \norm{\norm{\Box_{\ell} u}_{2}}_{\ell_{\ell}^{q}} = \norm{ \| \sum_{k \in \vee_{\ell}} \Box_{k}^{\al} \Box_{\ell} u\|_{2} }_{\ell_{\ell}^{q}}\\
			&= \norm{ \brk{\sum_{k \in \vee_{\ell}} \| \Box_{k}^{\al} u\|_{2}^{2}}^{1/2}}_{\ell_{\ell}^{q}}\\
			&\leq \norm{\norm{\|\Box_{k}^{\al} u\|_{2}}_{\ell_{k \in \vee_{\ell}}^{q}} \brk{\# \vee_{\ell}}^{1/2-1/q}}_{\ell_{\ell}^{q}}\\
			&\leq \norm{\inner{k}^{-\al d(1/2-1/q)/(1-\al)} \norm{\Box_{k}^{\al} u}_{2}}_{\ell_{k}^{q}} = \norm{u}_{M_{2,q}^{\al d(1/q-1/2)}}.
		\end{align*}
		\item[(b)] When $q\leq 1 \wedge p$, we have $\si(p,q) =0$. By the (quasi-)triangle inequality and the embedding of  $\ell^{r}$ spaces,  we have \begin{align*}
			\norm{u}_{M_{p,q}} &= \norm{\norm{\Box_{\ell} u}_{p}}_{\ell_{\ell}^{q}} = \norm{ \| \sum_{k \in \vee_{\ell}} \Box_{k}^{\al} \Box_{\ell} u\|_{p} }_{\ell_{\ell}^{q}}\\
			& \leq \norm{ \norm{ \| \Box_{k}^{\al }\Box_{\ell} u\|_{p}}_{\ell_{k \in \vee_{\ell}} ^{p\wedge 1}} }_{\ell_{\ell}^{q}} \\
			&\leq \norm{ \norm{\|\Box_{k}^{\al} u\|_{p}}_{\ell_{k \in \vee_{\ell}} ^{p\wedge 1}}} _{\ell_{\ell}^{q}} \lesssim \norm{u}_{\mpq^{0,\al}}.
		\end{align*}
		\item[(c)] When $(p, q) = (\infty, infty)$, we have $\si(p,q) = -d$. By the triangle inequality, we have 
		\begin{align*}
			\norm{u}_{M_{\infty,\infty}} &= \norm{\norm{\Box_{\ell} u}_{\infty}}_{\ell_{\ell}^{\infty}} = \norm{ \| \sum_{k \in \vee_{\ell}} \Box_{k}^{\al} \Box_{\ell} u\|_{\infty} }_{\ell_{\ell}^{\infty}}\\
			&\leq \norm{\sup_{k\in \vee_{\ell} }\norm{ \Box_{k}^{\al} u }_{\infty} \# \vee_{\ell} }_{\ell_{\ell}^{\infty}} \\
			&\leq \norm{\norm{\Box_{k}^{\al} u}_{\infty} \inner{k}^{-\al d/(1-\al)}}_{\ell_{k}^{\infty}} = \norm{u}_{M_{\infty,\infty}^{-\al d,\al}}.
		\end{align*}
		\item[(d)] When $p<1,q=\infty$, we have $\si(p,q) = -d/p$. By the quasi-triangular inequality, we have \begin{align*}
			\norm{u}_{M_{p,\infty}} &= \norm{\norm{\Box_{\ell} u}_{p}}_{\ell_{\ell}^{\infty}} = \norm{ \| \sum_{k \in \vee_{\ell}} \Box_{k}^{\al} \Box_{\ell} u\|_{p} }_{\ell_{\ell}^{\infty}}\\
			&\leq \norm{ \norm{ \| \Box_{k}^{\al} \Box_{\ell} u\|_{p}}_{\ell_{k \in \vee_{\ell}} ^{p}} }_{\ell_{\ell}^{\infty}}\\
			&\leq \norm{ \sup_{k\in \vee_{\ell} } \norm{\Box_{k}^{\al} u}_{p} \brk{\#\vee_{\ell}}^{1/p}} _{\ell_{\ell}^{\infty}} \\
			& \leq \norm{ \norm{\Box_{k}^{\al} u}_{p} \inner{k}^{-\al d/(p(1-\al)}}_{\ell_{k}^{\infty}} = \norm{u}_{M_{p,\infty}^{-\al d/p,\al}}.
		\end{align*}
		
	\end{description}
\end{proof}
		
\end{appendices}

\vskip .5cm
	\bibliographystyle{unsrt} 	
	\bibliography{smoothing}
\end{document}